\documentclass[a4paper]{amsart}
\usepackage[utf8]{inputenc}
\usepackage[english]{babel}

\usepackage{amscd,amssymb}
\usepackage{amsmath}
\usepackage{amsthm}
\usepackage{graphicx}
\usepackage{fancybox}
\usepackage{epic,eepic}
\usepackage{amstext}
\usepackage{mathtools}
\usepackage{esint}
\usepackage[square,numbers]{natbib}
\usepackage{booktabs}
\usepackage[hide links]{hyperref}
\usepackage{comment}
\usepackage{mathtools}
\usepackage{tikz,pgfplots}
\usepackage{subcaption}
\usepackage[noabbrev, capitalise, nameinlink]{cleveref}
\crefname{equation}{}{}
\crefname{assumption}{Assumption}{Assumptions}
\crefformat{equation}{\textup{#2(#1)#3}}

\newtheorem{theorem}{Theorem}[section]
\newtheorem{corollary}[theorem]{Corollary}
\newtheorem{lemma}[theorem]{Lemma}

\theoremstyle{definition}

\newtheorem{assumption}[theorem]{Assumption}
\theoremstyle{remark}
\newtheorem{remark}[theorem]{Remark}
\numberwithin{theorem}{section}
\numberwithin{equation}{section}
\numberwithin{figure}{section}

\def\with{\,:\,}
\def\dx{\,\mathrm{d}x}

\DeclareMathOperator*{\argmin}{arg\,min}

\newcommand{\jump}[1]{{[\![ #1 ]\!]}}

\newcommand{\setNodeDir}{\ensuremath{\setNode_\textup D}}

\renewcommand*{\vec}{\boldsymbol}

\newcommand*{\graph}{\ensuremath{\mathcal G}}
\newcommand*{\setEdge}{\ensuremath{\mathcal E}}
\newcommand*{\setNode}{\ensuremath{\mathcal N}}
\newcommand*{\edge}{\ensuremath{\mathfrak e}}
\newcommand*{\node}{\ensuremath{\mathfrak n}}
\newcommand{\halfpos}{n+\frac{1}{2}}
\newcommand{\halfneg}{n-\frac{1}{2}}
\newcommand{\Dt}{D_\tau}

\newcommand{\delx}{\partial_{\vec x}}

\newcommand{\urhn}{(\vec u_h^{n}, \vec r_h^{n})}
\newcommand{\urhnext}{(\vec u_h^{n+1}, \vec r_h^{n+1})}
\newcommand{\urhprev}{(\vec u_h^{n-1}, \vec r_h^{n-1})}
\newcommand{\urhtheta}{(\vec u_h^{n;\theta}, \vec r_h^{n;\theta})}
\newcommand{\urhgen}[1]{( \vec u_h^{#1}, \vec r_h^{#1})}

\newcommand{\Dturhgen}[1]{(\Dt \vec u_h^{#1}, \Dt \vec r_h^{#1})}
\newcommand{\urhdiff}[2]{(\vec u_h^{#1} - \vec u_h^{#2}, \vec r_h^{#1} - \vec r_h^{#2})}

\newcommand{\Con}[1][]{\ensuremath{C_{\textup{#1}}}}

\usepackage{todonotes}

\numberwithin{equation}{section}
\numberwithin{theorem}{section}

\AtBeginDocument{%
	\def\MR#1{}
}

\title[Numerical simulation of beam network models]{Numerical simulation of beam network models}
\author[M.~Görtz, M.~Hauck, A.~M{\aa}lqvist, A.~Rupp, L.~Swoboda]{Morgan~Görtz$^\text{\textsection}$, Moritz Hauck$^*$, Axel M{\aa}lqvist$^\ddagger$, Andreas Rupp$^\dagger$, Lucia~Swoboda$^\ddagger$}
\address{${}^*$ Institute for Applied and Numerical Mathematics, Karlsruhe Institute of Technology, Englerstr.~2, 76131 Karlsruhe, Germany}
\email{moritz.hauck@kit.edu}
\address{${}^{\dagger}$ Department of Mathematics, Faculty of Mathematics and Computer Science, Saarland University, Campus E1.1, 66123 Saarbrücken, Germany}
\email{andreas.rupp@uni-saarland.de}
\address{${}^{\ddagger}$ Department of Mathematical Sciences, Chalmers University of Technology \& University of Gothenburg, Chalmers Tvärgata 3, 412 96 Göteborg}
\email{axel@chalmers.se}
\email{lucias@chalmers.se}
\address{${}^{\text{\textsection}}$ Department of Computational Engineering and Design, Fraunhofer--Chalmers Centre, Chalmers Science Park, 412 88, Göteborg}
\email{morgan.gortz@fcc.chalmers.se}

\begin{document}

\begin{abstract}
	Network models are used as efficient representation of materials with complex, interconnected locally one-dimensional structures. They typically accurately capture the mechanical properties of a material, while substantially reducing computational cost by avoiding full three-dimensional resolution. Applications include the simulation of fiber-based materials, porous media, and biological systems such as vascular networks.
	This article focuses on two representative problems: a stationary formulation describing the elastic deformation of beam networks, and a time-dependent formulation modeling elastic wave propagation in such materials. We propose a two-level additive domain decomposition method to efficiently solve the linear system associated with the stationary problem, as well as the linear systems that arise at each time step of the time-dependent problem through implicit time discretization. We present a rigorous convergence analysis of the domain decomposition method when used as a preconditioner, quantifying the convergence rate with respect to network connectivity and heterogeneity. The efficiency and robustness of the proposed approach are demonstrated through numerical simulations of the mechanical properties of commercial-grade paperboard.
\end{abstract}


\maketitle

\section{Introduction}
Many problems in science and engineering, modeled by partial differential equations (PDEs), are posed in complex domains, which can lead to significant computational challenges. For example, when the computational domain is locally very slender, effectively one-dimensional, classical simulation techniques such as the finite element method (FEM) often require excessive spatial resolution and may produce nearly degenerate elements. As a computationally more tractable alternative to applying the FEM directly on such domains, one can employ dimension-reduced models of the underlying problem and discretize them.  For instance, one may represent the geometry of a material with a locally effectively one-dimensional structure by a graph, where edges correspond to the one-dimensional segments and nodes define their connections/joints. The mechanical behavior of these segments can then be described by one-dimensional  equations posed on the edges, which are suitably coupled at the network nodes. Such models, hereafter referred to as spatial network models, can substantially reduce computational complexity while preserving the essential characteristics of the full model.
 Applications where this modeling approach is feasible include porous media, in which the detailed pore geometry can be approximated by a network of nodes (pore cavities) and edges (throats) \cite{Blunt2013,Jivkov2013,Gjennestad2018,Huang2020}. Similarly, fiber-based materials such as paper can be represented as networks of one-dimensional beams corresponding to fibers, connected at nodes representing contact regions \cite{heyden2000network,Iliev2010,KMM20,Grtz2024}. This network-based modeling concept also extends to biological systems, such as vascular networks, where blood vessels can be represented as connected one-dimensional edges forming a spatial network \cite{Reichold2009,Mller2023,Fritz2022_1D0D3Dvascular}.

We emphasize that, despite using a locally dimension-reduced model, the linear systems resulting from discretization often remain challenging to solve efficiently due to very high condition numbers. This difficulty primarily stems from the multiscale character of the problem, reflected in strong material heterogeneity and complex network geometries. 
In the context of PDEs with multiscale coefficients, a setting closely related to our heterogeneous network problem, efficiently solving the resulting linear systems typically requires leveraging multiple discretization scales as well as parallelization. In what follows, we will highlight one prominent example of such methods, domain decomposition methods, and then briefly mention alternative approaches such as multiscale methods and (algebraic) multigrid methods.

Domain decomposition methods iteratively improve the current approximation by solving local error correction problems on subdomains, and update the current approximation with the computed correction. Each correction is obtained by solving a subdomain-restricted system matrix, using the residual restricted to that subdomain as the right-hand side. The resulting updates are combined either sequentially or simultaneously, corresponding to multiplicative and additive subspace correction methods, respectively. The subspaces in which the error corrections are solved are constructed by decomposing the computational domain into a collection of local subdomains. To prevent iteration counts from increasing with the number of subdomains, domain decomposition methods often incorporate a coarse level, i.e., the corresponding subspace correction equations are solved in a global low-dimensional, possibly problem-adapted, space. This allows information to be exchanged across distant subdomains, and the resulting method is referred to as a two-level scheme.
Two main families of domain decomposition methods are commonly distinguished: overlapping and non-overlapping methods.
In overlapping methods, subdomains partially overlap, enabling information exchange across the overlap regions. Examples of overlapping Schwarz methods using problem-adapted coarse spaces for improved treatment of heterogeneities include \cite{galvis2010domain-reduced,Nataf2011,efendiev2012robust,spillane2014abstract,gander2015analysis,Heinlein2019,heinlein2022fully,bastian2022multilevel,al2023efficient}. We also mention the works \cite{Strehlow2024,Alber2025}, which reformulate the Multiscale Spectral Generalized FEM, a multiscale method mentioned above, as a two-level overlapping domain decomposition method. 
Representative examples for non-overlapping domain decomposition methods include the Finite Element Tearing and Interconnecting (FETI) and Neumann–Neumann methods~\cite{Bjrstad2001,Spillane2013}, as well as their dual-primal and balancing variants, such as BDDC and FETI-DP~\cite{Mandel2007,Mandel2012,Klawonn2016,Pechstein2017,Kim2017}.

In the context of spatial networks, domain decomposition techniques have, for instance, been applied in \cite{Arioli2017} for a diffusion-type problem posed on a network, and in \cite{Leugering2017,Leugering2017b} for semilinear elliptic optimal control problems arising in gas network models. The latter works employ non-overlapping, one-level domain decomposition methods. Since no coarse space is included, the convergence rate generally deteriorates with an increased number of subdomains.
Two-level domain decomposition methods for spatial network models have been introduced in \cite{GoHeMa22,Grtz2024,Hauck2025}, where Timoshenko beam networks serve as model problems. In these works, coarse scales are artificially defined by a uniform Cartesian mesh overlaid on the network, with nodes assigned to mesh partitions accordingly.
The uniform convergence of the preconditioned iteration, originally established in the PDE setting, can be extended to the spatial network setting for discretization scales coarser than a certain critical threshold. In this regime, it is reasonable to assume specific homogeneity and connectivity properties of the underlying network. These properties enable the derivation of interpolation bounds that are central to the analysis, particularly for establishing the existence of a stable subspace decomposition required in the theory of additive Schwarz methods.

As an alternative to domain decomposition methods, also multiscale methods can be employed to obtain accurate and reliable approximations of multiscale problems. Their key idea is to incorporate microscopic coefficient information, for example through problem-adapted basis functions, by solving local fine-scale problems. This enables accurate approximations on coarse meshes that do not fully resolve the fine-scale structure of the coefficients. Typically, multiscale methods are one-shot approaches, obtaining an accurate solution within a single Galerkin approximation rather than iteratively improving it.
For diffusion-type PDEs, numerous multiscale methods construct problem-adapted basis functions, yielding optimal approximation orders under minimal assumptions on the coefficients. This comes at a moderate computational overhead compared to classical FEMs, typically due to enlarged basis function support or an increased number of basis functions per mesh entity. Prominent methods include the Generalized Multiscale FEM~\cite{EfeGH13,ChuEL18b}, the Multiscale Spectral Generalized FEM (MS-GFEM)~\cite{BabL11,Ma22}, Adaptive Local Bases~\cite{GraGS12}, the (Super-)Localized Orthogonal Decomposition (LOD) method~\cite{MalP14,HenP13,Malp20,HaPe21b,pumslod}, and Gamblets or operator-adapted wavelets~\cite{Owh17, OwhS19}; see also the review~\cite{AltHP21}.
Some concepts from PDE-based multiscale methods have been adapted to spatial network models. For example, \cite{EGHKM24,HaM22,Grtz2023,HMM23} apply the methodology of (Super-)LOD to spatial networks, defining coarse scales by superimposing a uniform Cartesian mesh. Homogenization-based upscaling methods have also been applied to heat-conducting networks \cite{Ewing2009,Iliev2010}, stochastic network problems \cite{Mansour2022}, traffic flow \cite{DellaRossa2010}, and porous media models \cite{CEPT12}.

Among available solvers for large linear systems from PDE discretizations, multigrid methods have proven particularly efficient. They recursively combine smoothing via simple iterative methods, such as Jacobi or Gauss–Seidel, which reduce high-frequency errors, with coarse-grid correction, which eliminates low-frequency errors \cite{Br77,Hackbusch1985}. For multiscale problems, classical geometric multigrid may fail to converge because standard smoothers and coarse spaces are not adapted to heterogeneous coefficients \cite{Alcouffe1981}. This challenge has motivated algebraic multigrid (AMG) methods, where coarse scales and prolongation operators are derived from the system matrix’s sparsity and connection strengths. Notable variants include Smoothed Aggregation AMG~\cite{Vank1996} and Energy-Minimizing AMG~\cite{Mandel1999,Wan1999}, which can handle highly heterogeneous coefficients \cite{XZ17}.

\subsection*{Outline}

In \cref{sec:timo}, we introduce the two model problems considered in this article: stationary and time-dependent elastic beam network problems. \cref{sec:disc} then presents their discretization using a classical high-order finite element method applied to network problems. A two-level domain decomposition preconditioner for the resulting, typically ill-conditioned, linear systems of equations is introduced in \cref{sec:precond}. The numerical experiments presented in \cref{sec:numexp} demonstrate the effectiveness of the proposed preconditioner for both model problems. Concluding remarks and potential avenues for future research are discussed in \cref{sec:conclusion}.

\section{Timoshenko beam network problems}\label{sec:timo}

The Timoshenko beam theory, originating from the work of Timoshenko~\cite{Timoshenko1921}, is widely used to describe the elastic deformation of beams.
It extends the classical Euler--Bernoulli beam theory by accounting for shear deformation. This broadens its applicability, making it suitable not only for long, slender beams but also for short, thick beams, where cross-sectional rotation cannot be neglected.
In this article, we consider two model problems:  
(i) a stationary problem describing the elastic deformation of a fiber network, and  
(ii) a time-dependent problem modeling the propagation of elastic waves in a fiber network. These model problems are subsequently introduced in two separate subsections.

\subsection{Stationary elasticity}
We represent the beam network as a graph $\graph = (\setNode, \setEdge)$, where $\setNode$ is the set of zero-dimensional nodes and $\setEdge$ the set of locally one-dimensional, straight edges.
The nodes correspond to the connection points between fibers, while the edges represent the fiber segments connecting two nodes.  
Throughout this work, we assume that the graph $\graph$ is connected.
The geometric domain of the network is abbreviated by
\begin{equation}
\label{eq:sigma}
    \Sigma \coloneqq \bigcup_{\edge \in \setEdge} \overline{\edge},
\end{equation}
where $\overline{\cdot}$ denotes the closure of a set.
Furthermore, for each edge $\edge \in \setEdge$, we denote by $\vec{i}\colon \Sigma \to \mathbb R^3$ the edge-wise constant function that assigns to an edge $\edge$ the unit vector pointing in its direction (sign fixed).
The elastic behavior of the beams is described by the two coefficient fields $\vec{B}, \vec{C} \colon \Sigma \to \mathbb{R}^{3\times 3}$, which are assumed to be 
uniformly positive definite and bounded. That is, there exist constants $0 < b_{\mathrm{min}} \leq b_{\mathrm{max}} < \infty$ and $0 < c_{\mathrm{min}} \leq c_{\mathrm{max}} < \infty$ such that, for almost every $\vec{x} \in \Sigma$,
\begin{align}\label{eq:boundCnm}
    \begin{split}
    b_\mathrm{min}|\vec \xi|^2 &\leq (\vec B(\vec{x}) \vec \xi,\vec \xi) \leq b_\mathrm{max}|\vec \xi|^2,
    \quad \forall \vec \xi \in \mathbb R^3,\\
    c_\mathrm{min}|\vec \xi|^2 &\leq (\vec C(\vec{x}) \vec \xi,\vec \xi) \leq c_\mathrm{max}|\vec \xi|^2,
    \quad \forall \vec \xi \in \mathbb R^3,
    \end{split}
\end{align}
where $(\cdot,\cdot)$ is the Euclidean inner product and $|\cdot| \coloneqq \sqrt{(\cdot,\cdot)}$ its induced norm.

The \emph{governing equations} for the elastic deformation of a beam, 
expressed in terms of the displacement and rotation 
$\vec{u}, \vec{r} \colon \Sigma \to \mathbb{R}^3$, take for any 
$\edge \in \setEdge$ the form
\begin{subequations}
\label{eq:stationary}
\begin{alignat}{2}
    - \delx \big( \vec{B} ( \delx \vec{u} + \vec{i} \times \vec{r} ) \big) 
    &= \vec{f},  &\quad& \text{in } \edge, \\ \label{homsysf}
    - \delx \big( \vec{C}\, \delx \vec{r} \big) 
    - \vec{i} \times (\vec{B} ( \delx \vec{u} + \vec{i} \times \vec{r} ))
    &= \vec{g},  &\quad& \text{in } \edge, 
\end{alignat}
where $\vec{f}, \vec{g} \colon \Sigma \to \mathbb{R}^3$ denote given 
distributed forces and moments, $\delx$ is the derivative in the 
direction of the edge $\edge$ parameterized with unit speed, and 
$\times$ denotes the standard cross product in $\mathbb{R}^3$.

To couple these equations on each edge, we impose continuity and balance conditions at the nodes. The problem is closed by Dirichlet boundary conditions at nodes in the non-empty set $\setNodeDir \subset \setNode$. The \emph{continuity conditions} require that displacements and rotations at a node are identical for all incident edges. Specifically, for each node $\node \in \setNode \setminus \setNodeDir$
and any two edges $\edge, \edge^\prime \in \setEdge$ incident to $\node$, we impose:
\begin{equation}\label{EQ:timo_cont_hybrid}
 \vec u|_\edge (\node) = \vec u|_{\edge^\prime}(\node), \qquad 
 \vec r|_\edge (\node) = \vec r|_{\edge^\prime}(\node).
\end{equation}
The \emph{balance conditions} enforce the equilibrium of forces and moments at non-Dirichlet nodes. For each $\node \in \setNode \setminus \setNodeDir$,
these conditions read
\begin{equation}\label{EQ:timo_balance}
  \jump{\vec{B} \left( \partial_{x} \vec{u} + \vec{i} \times \vec{r} \right) \nu}(\node) = 0, 
  \qquad  
  \jump{\vec{C} \, \partial_{\vec x} \vec{r} \, \nu}(\node) = 0,
\end{equation}
where the jump operator $\jump{\cdot}(\node)$ denotes the sum over all values attained at $\node$, 
and~$\nu$ is a function that assigns to the endpoints of each edge values in $\{+1,-1\}$. 
Specifically, $\nu$ takes the value $+1$ at the endpoint of $\edge$ if $\vec{i}|_\edge$ 
points outward from that node.
The \emph{Dirichlet boundary conditions} specify, for each node $\node \in \setNodeDir$, that
\begin{equation}\label{EQ:timo_dir}
\vec u(\node) = \vec u_\mathrm{D}(\node), \qquad
\vec r(\node) = \vec r_\mathrm{D}(\node),
\end{equation}
where $\vec u_\mathrm{D}, \vec r_\mathrm{D} \colon \setNodeDir \to \mathbb{R}^3$ are prescribed functions defining the Dirichlet data.
\end{subequations}

\subsubsection{Weak formulation}
To derive a weak formulation of problem \cref{eq:stationary}, we first introduce the function space for the displacement and rotation variables.  We define
\begin{equation}
 \vec V \coloneqq \left\{ \vec v \in \mathcal C^0(\Sigma)^3 \with \vec v|_\edge \in H^1(\edge)^3 ,\, \forall \edge \in \setEdge,\;\,  \vec v(\node) = \vec 0, \, \forall \node \in \setNodeDir \right\},
\end{equation}
where $\mathcal C^0(\Sigma)$ denotes the space of continuous functions defined on the network domain $\Sigma$.
A Friedrichs-type inequality holds on the network:
\begin{equation}
	\label{eq:friedrichs}
	\int_\Sigma |\vec v|^2 \, \mathrm{d}x
	\;\le\; C_\mathrm{F}^2
	\int_\Sigma |\partial_{\vec x} \vec v|^2 \,\mathrm{d}x,
	\quad \vec v \in \vec V,
\end{equation}
where the constant $C_\mathrm{F} > 0$ depends on the geometry of the considered network, specifically, the maximal intrinsic distance of a point $\vec x \in \Sigma$ to a Dirichlet node $\node \in \setNode_\mathrm{D}$. The space~$\vec V$ can be equipped with the norm
\begin{equation}
	\|\vec v\|_{\vec V}
	\coloneqq \left( \int_\Sigma |\partial_{\vec x} \vec v|^2 \,\mathrm{d}x\right)^{1/2},
	\quad \vec v \in \vec V,
\end{equation}
which is equivalent to the full $H^1(\Sigma)^3$-norm by virtue of \cref{eq:friedrichs}.
In the same manner, we equip the space~$\vec V \times \vec V$ with the norm
\begin{equation}
	\|(\vec v, \vec w) \|_{\vec V \times \vec V}
	\coloneqq \big( \|\vec v\|_{\vec V}^2 + \|\vec w\|_{\vec V}^2\big)^{1/2},
	\quad (\vec v, \vec w) \in \vec V \times \vec V.
\end{equation}

The bilinear and linear forms of the stationary problem read:
\begin{align}\label{eq:a}
a((\vec u,\vec r),(\vec \varphi,\vec \psi)) &\coloneqq \int_\Sigma (\vec B (\partial_{\vec x} \vec u \hspace{-.15ex}+\hspace{-.15ex} \vec i\times \vec r)) \cdot (\partial_{\vec x} \vec \varphi \hspace{-.15ex}+\hspace{-.15ex} \vec i \times \vec \psi)  +  (\vec C \partial_{\vec x} \vec r) \cdot \partial_{\vec x} \vec \psi \dx,
    \\
    F(\vec \varphi,\vec \psi) &\coloneqq \int_\Sigma \vec f \cdot \vec \varphi  + \vec g\cdot \vec \psi \dx.\label{eq:F}
\end{align}

Given source terms $\vec f, \vec g \in L^2(\Sigma)^3$, the weak formulation of the stationary beam network elasticity problem seeks the pair $(\vec{u},\vec{r}) \in \vec{V} \times \vec V$ such that
\begin{equation}\label{eq:weakstationary}
    a\big((\vec{u},\vec{r}),(\vec{\varphi},\vec{\psi})\big)
    = F(\vec{\varphi},\vec{\psi}),
    \quad \forall\, (\vec{\varphi},\vec{\psi}) \in \vec{V} \times \vec V.
\end{equation}

\subsubsection{Well-posedness}
To establish the well-posedness of this weak formulation, we begin by proving the coercivity and continuity of the bilinear form~$a$.
\begin{lemma}[Coercivity and boundedness of $a$]\label{lem:coercivity}
We have
\begin{align}
    a((\vec u,\vec r),(\vec u,\vec r)) 
    &\geq a_\mathrm{min} \,\|(\vec u, \vec r)\|_{\vec V \times \vec V}^2, \label{eq:coercivity}\\  
    a((\vec u,\vec r),(\vec \varphi,\vec \psi)) 
    &\leq a_\mathrm{max} \,\|(\vec u, \vec r)\|_{\vec V \times \vec V} \,\|(\vec \varphi, \vec \psi)\|_{\vec V \times \vec V},\label{eq:boundedness}
\end{align}
for all $(\vec u, \vec r), (\vec \varphi,\vec \psi)\in  \vec V\times \vec V$, where the constants $a_\mathrm{min}, a_\mathrm{max} > 0$ depend only on $b_\mathrm{min}, b_\mathrm{max}, c_\mathrm{min}, c_\mathrm{max}$ from \cref{eq:boundCnm} and the Friedrichs constant $C_\mathrm{F}$ from \cref{eq:friedrichs}.
\end{lemma}
\begin{proof}
We begin proving the coercivity of the bilinear form $a$. 
Let $(\vec u, \vec r) \in \vec V \times \vec V$ be arbitrary. 
Applying the weighted Young's inequality, we obtain, for any $\epsilon \in (0,1)$,
\begin{align*}
    &a((\vec u,\vec r),(\vec u,\vec r))\\
    &\; = \int_\Sigma |\vec B^{1/2} \partial_{\vec x} \vec u|^2 + 2( \vec B^{1/2}  \partial_{\vec x} \vec u)\cdot(\vec B^{1/2}(\vec i\times \vec r))  + |\vec B^{1/2}(\vec i \times \vec r)|^2 +  |\vec C^{1/2} \partial_{\vec x} \vec r|^2 \dx\\
    &\; \geq (1-\epsilon) \int_\Sigma |\vec B^{1/2} \partial_{\vec x} \vec u|^2 \dx + \int_\Sigma |\vec C^{1/2} \partial_{\vec x} \vec r|^2 \dx - \frac{1-\epsilon}{\epsilon}\int_\Sigma |\vec B^{1/2}(\vec i \times \vec r)|^2 \dx.
\end{align*}
Using the Friedrichs inequality from \cref{eq:friedrichs} and the uniform coefficient bounds \cref{eq:boundCnm}, we can estimate the last term on the right-hand side of the latter inequality as
\begin{align*}
    \int_\Sigma |\vec B^{1/2}(\vec i \times \vec r)|^2 \dx \leq b_\mathrm{max}\int_\Sigma|\vec r|^2 \dx \leq C \int_\Sigma |\vec C^{1/2} \partial_{\vec x} \vec r|^2 \dx,
\end{align*}
where we abbreviated $C\coloneqq b_\mathrm{max}c_\mathrm{min}^{-1}C_\mathrm{F}^2 >0$. Choosing any $\epsilon \in (\tfrac{C}{1+C},\,1)$ and again invoking the uniform coefficient bounds from \cref{eq:boundCnm} establishes estimate \cref{eq:coercivity}.

The boundedness of $a$, cf.~\cref{eq:boundedness}, can be established by combining the Friedrichs inequality from \cref{eq:friedrichs} with the uniform coefficient bounds in \cref{eq:boundCnm}.
\end{proof}

Introducing the $L^2$-norm for a tuple of $L^2(\Sigma)^3$-functions as
\begin{equation*}
    \|(\vec f, \vec g)\|_{L^2 \times L^2} \coloneqq \big( \|\vec f\|_{L^2(\Sigma)^3}^2 + \|\vec g\|_{L^2(\Sigma)^3}^2 \big)^{1/2}, \quad (\vec f, \vec g) \in L^2(\Sigma)^3 \times L^2(\Sigma)^3,
\end{equation*}
we obtain the following well-posedness result.

\begin{theorem}[Well-posedness of stationary elasticity]
The weak formulation of the beam network problem \cref{eq:weakstationary} is well-posed. Specifically, there exist a unique solution $(\vec u, \vec r) \in \vec V \times \vec V$ satisfying \cref{eq:weakstationary} and the following stability estimate holds:
    \begin{equation}
    \label{eq:stability}
        \|(\vec{u},\vec{r})\|_{\vec V\times \vec V} \leq a_\mathrm{min}^{-1}C_\mathrm{F}\|(\vec{f},\vec{g})\|_{L^2\times L^2}.
    \end{equation}
    \end{theorem}
\begin{proof}
Let us first show the boundedness of the linear form $F$. Using Friedrichs inequality from \cref{eq:friedrichs}, we obtain the following bound:
\begin{equation}
\label{eq:rhsboundedness}
    \sup_{(\vec \varphi, \vec \psi) \in \vec V \times \vec V} \frac{F(\vec \varphi, \vec \psi)}{\|(\vec \varphi, \vec \psi)\|_{\vec V \times \vec V}} \leq C_\mathrm{F} \|(\vec f, \vec g)\|_{L^2 \times L^2}.
\end{equation}
Well-posedness then follows directly from the Riesz representation theorem, and the stability estimate \cref{eq:stability} is obtained by combining the coercivity estimate \cref{eq:coercivity}, the fact that $(\vec u, \vec r)$ solves \cref{eq:weakstationary}, and \cref{eq:rhsboundedness}.
\end{proof}

\subsection{Elastic wave propagation}
For the elastic wave propagation problem, in addition to the coefficients 
$\vec B$ and $\vec C$ from \cref{eq:boundCnm}, we introduce one further 
coefficient field $\vec D \colon \Sigma \to \mathbb{R}^{3\times 3}$. 
We assume that $\vec D$ is uniformly positive definite and bounded, 
i.e., there exist constants $0 < d_\mathrm{min} \le d_\mathrm{max} < \infty$ 
such that, for almost every $\vec x \in \Sigma$,
\begin{equation}\label{eq:boundD}
    d_\mathrm{min} |\vec \xi|^2 
    \le (\vec D(\vec x)\vec \xi, \vec \xi) 
    \le d_\mathrm{max} |\vec \xi|^2, 
    \quad \forall \vec \xi \in \mathbb{R}^3.
\end{equation}

We suppose that the elastic wave propagation problem is posed on the time interval $(0,T)$ with end time $T>0$. The \emph{governing equations} for the elastic wave propagation in a beam, expressed in terms of the (now time-dependent) displacement and rotation $\vec{u}, \vec{r}\colon (0,T) \times \Sigma \to \mathbb R^3$, then take for any $\edge \in \setEdge$ the form
\begin{subequations}\label{eq:elasticwave}
\begin{alignat}{2}
    c \ddot{\vec{u}} 
    - \delx \vec{B} \big( \delx \vec{u} + \vec{i} \times \vec{r} \big) 
    &= \vec{f}, &\quad& \text{in } (0,T)\times\edge, \label{homsysf} \\
    \vec{D} \ddot{\vec{r}} 
    - \delx ( \vec{C} \delx \vec{r} ) 
    - \vec{i} \times \vec{B} ( \delx \vec{u} + \vec{i} \times \vec{r} ) 
    &= \vec{g}, &\quad& \text{in } (0,T)\times\edge, \label{homsysg}
\end{alignat}
    where $\vec f, \vec g \colon (0,T) \times \Sigma \to \mathbb R^3$ are given distributed forces and moments, $c>0$ is a constant
    and $\ddot{\vec u}$ denotes the second-order time derivative of $\vec u$; the same notation applies for $\vec r$. In addition, analogously to the stationary problem, for almost all $t \in (0,T)$, we impose the coupling conditions \cref{EQ:timo_balance,EQ:timo_cont_hybrid}, and the Dirichlet boundary conditions \cref{EQ:timo_dir}. 
The system is supplemented with the initial conditions
\begin{equation}
\label{ic:weakwave}
    \vec u(0,\cdot) = \vec u_0, \quad 
    \dot{\vec u}(0,\cdot) = \vec v_0, \quad 
    \vec r(0,\cdot) = \vec r_0, \quad 
    \dot{\vec r}(0,\cdot) = \vec p_0.
\end{equation}
\end{subequations}

\subsubsection{Weak formulation}
To derive a weak formulation of the elastic wave propagation problem, 
we first introduce some additional notation for time-dependent function spaces. 
We identify a function 
$\vec v \colon (0,T) \times \Sigma \to \mathbb{R}^3$ 
with the corresponding mapping 
$\vec v \colon (0,T) \to \vec X$, 
where $\vec X$ is a suitable Banach space, and employ the standard definition and notation for the Bochner spaces $L^p(0,T;\vec X)$,  $1 \leq p \leq \infty$, 
see, e.g.,~\cite{Evans2010}. 
Moreover, we denote by $\vec V^*$ the dual of $\vec V$, 
and by $\langle \cdot,\cdot \rangle_{\vec V^* \times \vec V}$ 
the canonical duality pairing. 
For $\vec v, \vec q \in \vec V^*$ and 
$\vec \varphi, \vec \psi \in \vec V$, we define the bilinear form
\begin{equation}
\label{eq:m}
    m\big((\vec v,\vec q),(\vec \varphi, \vec \psi)\big) 
    \coloneqq \langle c\, \vec v, \vec \varphi \rangle_{\vec V^* \times \vec V} 
    + \langle \vec D\, \vec q, \vec \psi \rangle_{\vec V^* \times \vec V},
\end{equation}
where the constant $c$ was introduced above in \cref{homsysf}. Furthermore, we define the constants
$m_\mathrm{min} \coloneqq \min(c, d_\mathrm{min})$ and $m_\mathrm{max} \coloneqq \max(c, d_\mathrm{max})$. Note that \( \vec{V} \subset L^2(\Sigma)^3 \subset \vec{V}^\ast \) forms a Gelfand triple, and the duality pairing coincides with the \(L^2\)-inner product whenever both functions lie in \(L^2(\Sigma)^3\).

Given source terms $\vec f, \vec g \in L^2(0,T;L^2(\Sigma)^3)$,
the weak formulation of the elastic wave propagation problem 
seeks the pair  
$(\vec u, \vec r) \in L^2(0,T;\vec V \times \vec V)$ with first-order time derivatives 
$(\dot{\vec u}, \dot{\vec r}) \in L^2(0,T;L^2(\Sigma)^3 \times L^2(\Sigma)^3)$ and second-order 
time derivatives $(\ddot{\vec u}, \ddot{\vec r}) \in L^2(0,T;\vec V^* \times \vec V^*)$, 
such that, for almost all $t \in (0,T)$,
\begin{subequations}
\label{eq:weakwave}
\begin{equation}
    m\big((\ddot{\vec u}, \ddot{\vec r}), (\vec \varphi, \vec \psi)\big) 
    + a\big((\vec u, \vec r), (\vec \varphi, \vec \psi)\big) 
    = F(\vec \varphi, \vec \psi), 
    \quad \forall (\vec \varphi, \vec \psi) \in \vec V \times \vec V,
\end{equation}
with the initial conditions \cref{ic:weakwave}, now expressed as
\begin{equation}
\label{ic:weakwave2}
    (\vec u,\vec r)(0) = (\vec u_0, \vec r_0), \quad 
    (\dot{\vec u},\dot{\vec r})(0) = (\vec v_0,\vec p_0),
\end{equation}
where $\vec u_0, \vec r_0 \in \vec V$ and $\vec v_0, \vec p_0 \in L^2(\Sigma)^3$.
\end{subequations}

\subsubsection{Well-posedness}
The well-posedness result stated in the following theorem relies on energy techniques. 
In the present context, the energy is defined for a pair of functions 
$(\vec{v}, \vec{q}) \in L^\infty(0,T;\vec V \times \vec V)$ with first-order time derivatives 
$(\dot{\vec{v}}, \dot{\vec{q}}) \in L^\infty(0,T;L^2(\Sigma)^3\times L^2(\Sigma)^3)$, for almost every $t \in (0,T)$, by
\begin{align*}
    \mathcal{E} (\vec{v}, \vec{q})(t) \coloneqq 
    \frac{1}{2} m\big((\dot{\vec{v}}(t),\dot{\vec{q}}(t)), (\dot{\vec{v}}(t),\dot{\vec{q}}(t))\big)
    + \frac{1}{2} a\big((\vec{v}(t), \vec{q}(t)), (\vec{v}(t), \vec{q}(t))\big).
\end{align*}

\begin{theorem}[Well-posedness of elastic wave propagation]
The weak formulation \cref{eq:weakwave} is well-posed. In particular, there exists a unique solution $(\vec u, \vec r) \in L^\infty(0,T; \vec V \times \vec V)$ with $(\dot{\vec u}, \dot{\vec r}) \in L^\infty(0,T; L^2(\Sigma)^3 \times L^2(\Sigma)^3)$ and $(\ddot{\vec u}, \ddot{\vec r}) \in L^2(0,T; \vec V^* \times \vec V^*)$, satisfying \cref{eq:weakwave}. Moreover, the following stability estimate holds:
    \begin{align}
    &\| (\vec u , \vec r) \|_{L^\infty(0,T; \vec V \times \vec V)} 
    + \| (\dot{\vec u}, \dot{\vec r}) \|_{L^\infty(0,T; L^2 \times L^2)}
    \nonumber\\
    &\leq C (\| (\vec{v}_{0}, \vec{p}_{0}) \|_{L^2 \times L^2}
     + \| (\vec{u}_{0}, \vec{r}_{0}) \|_{\vec V \times \vec V}
     + \| (\vec{f}, \vec{g}) \|_{L^2(0,T; L^2 \times L^2)}),
    \end{align}
where the constant $C$ depends only on the end time $T$ and the constants $a_\mathrm{min},$ $a_\mathrm{max},$ $m_\mathrm{min},$ $m_\mathrm{max}$.
\end{theorem}
\begin{proof}
    To prove the existence of a solution, we consider a Galerkin approximation in the space $L^2(0,T; \vec V_h \times \vec V_h)$,
    where $\vec V_h$ denotes a conforming finite element subspace of $\vec V$, and prove the existence of the solution for the discrete problem. We choose a finite-dimensional space~$\vec V_h$ of dimension $N \in \mathbb N$ such that 
    \begin{align*}
        \lim_{h \to 0} \inf_{\vec v_h \in \vec V_h} \| \vec v - \vec v_h\|_{\vec V} = 0, \quad \forall \vec v \in \vec V.
    \end{align*}
This yields the finite-dimensional problem: find $(\vec{u}_h, \vec{r}_h) \in L^2(0,T;\vec V_h \times \vec V_h)$, such that, for almost all $t \in (0,T)$ and all $(\vec{\varphi}_h, \vec{\psi}_h) \in \vec V_h \times \vec V_h$,
\begin{align}
    m\big((\ddot{\vec{u}}_h(t), \ddot{\vec{r}}_h(t)),(\vec{\varphi}_h, \vec{\psi}_h)\big)
    + a\big((\vec{u}_h(t), \vec{r}_h(t)),(\vec{\varphi}_h, \vec{\psi}_h)\big) = F(\vec \varphi_h,\vec \psi_h), \label{weakd}
\end{align}
with the initial conditions
\begin{align}
    (\vec{u}_h, \vec{r}_h)(0) = (\vec{u}_{0h}, \vec{r}_{0h}), \qquad (\dot{\vec{u}}_h, \dot{\vec{r}}_h)(0) = (\vec{v}_{0h},\vec{p}_{0h}), \label{iconddiscrete}
\end{align}
where $\vec{u}_{0h}$, $ \vec{r}_{0h}$, $\vec{v}_{0h}$, and $\vec{p}_{0h}$ are defined as
\begin{align*}
    \vec{u}_{0h} \coloneqq \argmin_{\vec{w}_h \in \vec V_h} \| \vec{u}_{0} - \vec{w}_h \|_{L^2(\Sigma)^3}, \qquad 
    \vec{r}_{0h} \coloneqq \argmin_{\vec{w}_h \in \vec V_h} \| \vec{r}_{0} - \vec{w}_h \|_{L^2(\Sigma)^3},\\
    \vec{v}_{0h} \coloneqq \argmin_{\vec{w}_h \in \vec V_h} \| \vec{v}_{0} - \vec{w}_h \|_{L^2(\Sigma)^3}, \qquad 
    \vec{p}_{0h} \coloneqq \argmin_{\vec{w}_h \in \vec V_h} \| \vec{p}_{0} - \vec{w}_h \|_{L^2(\Sigma)^3}.
\end{align*}

Problem \cref{weakd} can be equivalently written in matrix-vector form as
\begin{align*}
    \vec{M}_h \ddot{\vec{\xi}}(t)
    + \vec{A}_h \vec{\xi}(t)  = \vec{F}_h(t),
\end{align*}
where $\vec{\xi}(0)$ is the coefficient vector of $(\vec u_{0h}, \vec r_{0h})$ and $\dot{\vec{\xi}}(0)$ is the coefficient vector of $(\vec v_{0h}, \vec p_{0h})$. Since $\vec{M}_h$ is invertible, this system admits a unique solution with absolutely continuous time derivative. The energy estimates of the discrete solution follow from the conservation of energy and an application of Gronwall’s inequality in differential form, which yields the estimate
\begin{align*}
    \mathcal{E}({\vec{u}}_h, {\vec{r}}_h)(t) \leq
    e^{m_\mathrm{min}^{-1} t} \left( \mathcal{E}({\vec{u}}_h, {\vec{r}}_h)(0) + \frac{1}{2} \int_0^t \|(\vec{f}(s), \vec{g}(s)) \|_{L^2 \times L^2}^2 \,ds \right).
\end{align*}
We obtain
\begin{align*}
    \frac{1}{2} \mathrm{min}(a_\mathrm{min}, m_\mathrm{min}) \left( \| (\dot{\vec{u}}_{h}, \dot{\vec{r}}_{h}) \|_{L^2 \times L^2}^2 + \| ( {\vec{u}}_{h}, {\vec{r}}_{h}) \|_{\vec V \times \vec V}^2 \right) \leq \mathcal{E}({\vec{u}}_h, {\vec{r}}_h)(t),
\end{align*}
and
\begin{align*}
    \max_{t \in [0,T]} \mathcal{E}({\vec{u}}_h, {\vec{r}}_h)(t) \leq \alpha
    (\| (\vec{v}_{0}, \vec{p}_{0}) \|_{L^2 \times L^2}^2
     + \| (\vec{u}_{0}, \vec{r}_{0}) \|_{\vec V \times \vec V}^2
     + \| (\vec{f}, \vec{g}) \|_{L^2(0,T; L^2 \times L^2)}^2),
\end{align*}
where $\alpha \coloneqq \frac{1}{2} e^{m_\mathrm{min}^{-1} T} \mathrm{max}(a_\mathrm{max}, m_\mathrm{max}, 1).$
Moreover, we conclude the boundedness of \((\ddot{\vec{u}}_{h},\, \ddot{\vec{r}}_{h})\).
Namely, for $(\vec \varphi_h, \vec \psi_h)$ denoting the $L^2$-projection of $(\vec \varphi, \vec \psi)$ onto $\vec V_h \times \vec V_h$, we obtain the estimate
\begin{align*}
    \| (\ddot{\vec{u}}_h(t), \ddot{\vec{r}}_h(t))\|_{\vec V^* \times \vec V^*}
    &\lesssim ||(\vec f(t), \vec g(t))||_{L^2 \times L^2} + ||(\vec u_h(t), \vec r_h(t))||_{\vec V \times \vec V}.
\end{align*}
Here, we have used Friedrichs' inequality and the $H^1$-stability of the $L^2$-projection in one dimension.
Since the families $(\vec{u}_h, \vec{r}_h)_{h>0}$, $(\dot{\vec{u}}_h, \dot{\vec{r}}_h)_{h>0}$ and $(\ddot{\vec{u}}_h, \ddot{\vec{r}}_h)_{h>0}$ are bounded in $L^2(0,T;\vec V \times \vec V)$, $L^2(0,T;L^2 \times L^2)$, and $L^2(0,T;\vec V^* \times \vec V^*)$, respectively, there exists a subsequence that converges weakly in these spaces.
Assuming nested meshes, it follows by a standard limit argument that the weak limit is a solution of problem~\cref{eq:weakwave}. The uniqueness of the solution follows since, for $\vec{f} = \vec{g} = 0$ and zero initial data, the only weak solution is $\vec{u} = \vec{r} = 0$.
\end{proof}

A direct consequence of this proof is the following energy conservation property.
\begin{corollary}[Energy conservation]
\label{cor:energy-cons}
For the solution pair \((\vec u, \vec r)\) to problem~\cref{eq:weakwave}, 
the following holds for almost every \(t \in (0,T)\):
    \begin{equation*}
        \frac{\mathrm{d}}{\mathrm{d}t} \mathcal E(\vec u,\vec r)(t) = \int_\Sigma \vec f(t) \cdot \dot{\vec u}(t)  + \vec g(t) \cdot \dot{\vec r}(t) \dx.
    \end{equation*}
    Therefore, for $\vec f = \vec g = 0$, the initial energy of the system is preserved over time.
\end{corollary}

\section{Discretization}\label{sec:disc}

In this section, we first introduce a high-order finite element discretization of the stationary problem. For the elastic wave propagation problem, we then combine this space discretization with an appropriate time discretization.

\subsection{Stationary elasticity}

As usual, a finite element discretization is based on a mesh $\mathcal T_h$ of the domain.  
In our network setting, the domain is $\Sigma$ defined in \cref{eq:sigma}, and the mesh $\mathcal T_h$ is obtained by subdividing each fiber into elements.  
The mesh parameter $h > 0$ is the maximum element length in $\mathcal T_h$.  
For a fixed polynomial degree $p \ge 1$, the  discrete approximation space associated with $\mathcal T_h$ is given by
\begin{align}
    \vec{V}_h \coloneqq \big\{ \vec v \in \vec V \; \with \; \vec v|_K \in (\mathbb P^p(K))^3, \ \forall\, K \in \mathcal T_h \big\}, \label{DGspace}
\end{align}
where $\mathbb P^p(K)$ denotes the space of univariate polynomials on $K$ of degree at most~$p$.  

The finite element approximation of the network elasticity problem \cref{eq:weakstationary} then seeks a discrete pair $(\vec u_h,\vec r_h) \in \vec V_h \times \vec V_h$ satisfying
\begin{equation}\label{eq:weakstationarydisc}
    a\big((\vec{u}_h,\vec{r}_h),(\vec{\varphi}_h,\vec{\psi}_h)\big)
    = F(\vec{\varphi}_h,\vec{\psi}_h),
    \quad \forall\, (\vec{\varphi}_h,\vec{\psi}_h) \in \vec{V}_h \times \vec V_h.
\end{equation}
Note that, by the conformity of $\vec V_h$, the well-posedness of \cref{eq:weakstationarydisc} follows directly from the continuous case, and also the stability estimate \cref{eq:stability} remains valid in the discrete setting with the same constant.

In the convergence analysis of the finite element method, the best-approximation error is estimated by introducing a suitable interpolation operator 
$\mathcal I_h \colon \vec V \to \vec V_h$, satisfying the following approximation error estimate:
\begin{equation}\label{eq:Ih}
    \|\vec v - \mathcal I_h \vec v\|_{\vec V} \leq C_{\mathcal I} h^{p} |\vec v|_{H^{p+1}(\mathcal T_h)^3}, 
    \quad \vec v \in \vec V \cap H^{p+1}(\mathcal T_h)^3,
\end{equation}
where $C_{\mathcal I} > 0$ is a constant independent of $h$ and $|\cdot|_{H^k(\mathcal T_h)^3}$ denotes the broken $H^k$-seminorm; see, e.g., \cite[Cor.~1.110]{Ern2004}. Introducing the notation,
\begin{equation*}
    |(\vec u,\vec r)|_{H^k\times H^k} \coloneqq 
    \big(|\vec u|_{H^k(\mathcal T_h)^3}^2 + |\vec r|_{H^k(\mathcal T_h)^3}^2\big)^{1/2},
\end{equation*}
we obtain the following convergence theorem.

\begin{theorem}[Convergence for stationary elasticity] \label{thm:conv_stat}
Assume that the solution pair satisfies the piecewise regularity requirement $(\vec u,\vec r) \in H^{p+1}(\mathcal T_h)^3 \times H^{p+1}(\mathcal T_h)^3$. Then, we have the following error estimate:
    \begin{align}
        \|(\vec u-\vec u_h,\vec r-\vec r_h)\|_{\vec V\times \Vec V} \leq a_\mathrm{min}^{-1}a_\mathrm{max} C_{\mathcal I} h^p|(\vec u,\vec r)|_{H^{p+1}\times H^{p+1}}.
    \end{align}
\end{theorem}
\begin{proof}
The convergence analysis follows the classical approach. First, Céa's lemma bounds the approximation error of the finite element solution, cf.~\cref{eq:weakstationarydisc}, by the best-approximation error in the space $\vec V_h \times \vec V_h$. Applying the interpolation error estimate \cref{eq:Ih} then yields the desired error bound.
\end{proof}

\subsection{Elastic wave propagation}

Next, we discuss the discretization of the elastic wave propagation problem \cref{eq:elasticwave}. We begin by performing a spatial discretization, followed by a temporal discretization. Applying the same spatial discretization as for the stationary problem above, we obtain the following semi-discrete formulation: Seek $(\vec{u}_h,\vec{r}_h) \in H^2(0,T;\vec{V}_h \times \vec{V}_h)$ such that, for almost all $t \in (0,T)$,
\begin{equation}\label{eq:weakwavesemidisc}
    m\big((\ddot{\vec{u}}_h,\ddot{\vec{r}}_h),(\vec{\varphi}_h,\vec{\psi}_h)\big)+
    a\big((\vec{u}_h,\vec{r}_h),(\vec{\varphi}_h,\vec{\psi}_h)\big)
    = F(\vec{\varphi}_h,\vec{\psi}_h),
\end{equation}
for all $(\vec{\varphi}_h,\vec{\psi}_h) \in \vec{V}_h \times \vec V_h$.

To introduce the temporal  discretization, we first fix some notation. Let $\{\vec v^n\}_{n=0}^N$ be a sequence of functions with $\vec v^n \in \vec V$ for $n \in \{0, \dots, N\}$. We define the average of two consecutive functions by
\begin{equation*}
	\vec v^{n+\frac{1}{2}}\coloneqq\tfrac{\vec v^{n+1}+ \vec v^{n}}{2},
\end{equation*}
for $n \in \{0,\dots,N-1\}$. 
For such a sequence and the uniform time step $\tau=\frac{T}{N}$, we define the first- and second-order discrete time derivatives by
\begin{equation*}
	\begin{aligned}
		\Dt \vec v^{\halfpos} \coloneqq\tfrac{ \vec v^{n+1}- \vec v^n}{\tau},\qquad
		\Dt^2 \vec v^{n} \coloneqq\tfrac{\vec v^{n+1}-2 \vec v^n+ \vec v^{n-1}}{\tau^2}.
	\end{aligned}
\end{equation*}

In the following, the initial conditions are given by
\begin{equation} \label{eq:intial_condition}
    (\vec u_h^0, \vec r_h^0) \coloneqq P_m(\vec u_{0}, \vec r_{0}),
\end{equation}
where $P_m$ denotes the $m$-orthogonal projection onto the finite element space. Furthermore, the fictitious pair $(\vec u_h^{-1}, \vec r_h^{-1})$ is defined implicitly through
\begin{equation} \label{eq:central_diff_tau}
    \frac{( \vec u_h^1, \vec r_h^1) - ( \vec u_h^{-1}, \vec r_h^{-1})}{2 \tau} \coloneqq P_m(\vec v_{0}, \vec p_{0})
\end{equation}
for $\theta \neq \tfrac{1}{12}$, and by
\begin{equation} \label{eq:central_diff_tau12}
    \frac{( \vec u_h^1, \vec r_h^1) - ( \vec u_h^{-1}, \vec r_h^{-1})}{2 \tau} \coloneqq P_m(\vec v_{0}, \vec p_{0}) + \frac{\tau^2}{6} P_m(\partial_t^3 \vec u(0), \partial_t^3 \vec r(0))
\end{equation}
for $\theta = \tfrac{1}{12}$.
We employ the $\theta$-scheme for the time discretization of the semi-discrete wave propagation problem using a uniform time step; see, e.g., \cite{Jo03, Maier25}.
This leads to the following discrete formulation. Find
$\{(\vec u_h^n, \vec r_h^n)\}_{n=0}^{N}$ with $(\vec u_h^n, \vec r_h^n) \in \vec V_h \times \vec V_h$, such that, for $1 \le n \le N$, 
\begin{equation}\label{eq:th}
 m\big( (\Dt^2 \vec u_h^{n}, \Dt^2 \vec r_h^{n}), (\vec \varphi_h, \vec \psi_h)\big)  
 + a\big( \urhtheta, (\vec \varphi_h, \vec \psi_h) \big) 
 = F^{n;\theta}(\vec \varphi_h, \vec \psi_h), 
\end{equation}
for all $(\vec \varphi_h, \vec \psi_h) \in \vec V_h \times \vec V_h$ and
appropriately chosen
initial data, where $\theta\in[0,\frac{1}{2}]$ is fixed. The weighted $\theta$-difference is defined by
\begin{equation}\label{eq:theta_difference}
\urhtheta \coloneqq
 \theta \urhnext
+ (1 - 2\theta)\urhn
+ \theta \urhprev.
\end{equation}
Further, the right-hand side functional $F^{n;\theta}$ is given by
\begin{equation}\label{eq:Fnt}
    F^{n;\theta}(\vec \varphi_h, \vec \psi_h) \coloneqq
    \int_\Sigma \vec f^{n;\theta} \cdot \vec \varphi_h  
    + \vec g^{n;\theta} \cdot \vec \psi_h \, \mathrm{d}x,
\end{equation}
where $\vec f^{n;\theta}$ and $\vec g^{n;\theta}$ are defined analogously to~\eqref{eq:theta_difference} with 
$\vec f^n = \vec f(n\tau)$ and $\vec g^n = \vec g(n\tau)$, provided that $\vec f$ and $\vec g$ can be evaluated pointwise.  
Note that this time discretization generalizes classical schemes: for $\theta=0$ it reduces to the leapfrog method, while $\theta=\tfrac{1}{4}$ yields the Crank--Nicolson scheme.

\subsubsection{Stability}
The sequence of discrete solution pairs satisfies, analogously to the time-continuous case in \cref{cor:energy-cons}, certain energy conservation or stability properties, depending on the specific choice of~$\theta$.
For the $\theta$-method~\eqref{eq:th} as temporal discretization, we define the discrete energy
\begin{equation}\label{eq:nrg_definition}
\begin{aligned}
\mathcal{E}^{\halfpos}
 & \coloneqq \tfrac{1}{2}\Big[ 
  m\big( \Dturhgen{\halfpos}, \Dturhgen{\halfpos}\big) \\
  &+ a\big(\urhgen{\halfpos}, \urhgen{\halfpos}\big) \\
  &+ \tau^2\!\left(\theta-\tfrac{1}{4}\right) 
    a\big( \Dturhgen{\halfpos}, \Dturhgen{\halfpos}\big)
\Big],
\end{aligned}
\end{equation}
where the last term is a numerical correction appearing for $\theta \neq \frac{1}{4}$. With regard to the discrete energy conservation, it suffices to show that the discrete energy of the $\theta$-method is non-negative at all times to obtain stability.

\begin{theorem}[Stability]\label{thm:stability}
For~$\tfrac{1}{4}\leq\theta\leq\tfrac{1}{2}$, the $\theta$-method \eqref{eq:th} is unconditionally stable. For~$0\leq\theta<\tfrac{1}{4}$, the $\theta$-method \eqref{eq:th} is stable if the CFL condition
\begin{equation}\label{eq:CFL}
	\tau \leq \frac{1-\delta}{\Con[CFL]} h_\mathrm{min},\quad \Con[CFL]= \sqrt{\frac{a_{\mathrm{max}}}{m_\mathrm{min}}\left( \frac{1}{4} - \theta\right)} 
 C_\mathrm{inv}
\end{equation}
holds for some $\delta \in (0,1)$, where $C_\mathrm{inv}$ is the constant from the inverse inequality
\begin{align*}
    \| (\vec v_h, \vec w_h) \|_{\vec V \times \vec V} \leq C_\mathrm{inv} h_\mathrm{min}^{-1} \| (\vec v_h, \vec w_h) \|_{L^2 \times L^2}.
\end{align*}
In both cases, there exists a constant $\Con[s]>0$
depending on $\delta, a_{\textrm{min}}, a_{\textrm{max}}, m_{\textrm{min}}$ and $m_{\textrm{max}}$ such that, for all $n \in \{0,\dots,N-1\}$, 
\begin{equation}\label{eq:stability_estimate}
\begin{aligned}
\big\| &\Dturhgen{\halfpos}\big\|_{L^2 \times L^2}+\big\|\urhgen{\halfpos}\big\|_{\vec V \times \vec V}\\
&\leq  \Con[s]\Big(\big\| \Dturhgen{\frac{1}{2}}\big\|_{L^2 \times L^2}+\sqrt{\|\urhgen{1}\|_{\vec V \times \vec V}\|\urhgen{0}\|_{\vec V \times \vec V}}\\
&\qquad\qquad+\tau\sqrt{\theta}\big\| \Dturhgen{\frac{1}{2}}\big\|_{\vec V \times \vec V} +\sum_{k=1}^n\tau\big\|( \vec f^{k;\theta}, \vec g^{k;\theta} )\big\|_{L^2 \times L^2}\Big).
\end{aligned}
\end{equation}
\end{theorem}

\begin{proof}
    For $\theta \ge \frac{1}{4}$, it is immediate that the method is unconditionally stable, since the last term in \cref{eq:nrg_definition} is nonnegative. For $\theta < \frac{1}{4}$, the CFL condition yields
    \begin{equation}
        \mathcal{E}^{\halfpos} \geq \frac{1}{2} a\big(\urhgen{\halfpos}, \urhgen{\halfpos}\big) + \frac{1}{2} c_\delta \|\urhgen{\halfpos}\|_{L^2 \times L^2}^2 \geq 0,
    \end{equation}
    for some $\delta >0$, where $c_\delta \coloneqq m_\mathrm{min} (1-(1-\delta)^2)$. Furthermore, we obtain
    \begin{equation}
        \big\| \Dturhgen{\halfpos}\big\|_{L^2 \times L^2}+\big\|\urhgen{\halfpos}\big\|_{\vec V \times \vec V}
        \leq 2 \alpha \sqrt{\mathcal{E}^{\halfpos}},
    \end{equation}
    where $\alpha \coloneqq \sqrt{\max(a_{\textrm{min}}^{-1}, c_\delta^{-1})}.$ The inequality
    \begin{equation}
        \sqrt{\mathcal{E}^{\halfpos}} - \sqrt{\mathcal{E}^{\halfneg}}
        \leq \tau \alpha \big\|( \vec f^{n;\theta}, \vec g^{n;\theta} )\big\|_{L^2 \times L^2}
    \end{equation}
    and a bound on the initial energy (that depends on $a_{\textrm{max}}$, $m_{\textrm{max}}$, $\tau$ and $\theta$) completes the proof.
\end{proof}

\begin{corollary}[Discrete energy conservation]\label{cor:nrg}
The discrete energy
$\mathcal{E}^{\halfpos}$
satisfies
\begin{equation}\label{eq:nrg_estimate}
 F^{n;\theta}(\vec u_h^{\,n+1} - \vec u_h^{\,n-1}, \, \vec r_h^{\,n+1} - \vec r_h^{\,n-1})
= 2\big(\mathcal{E}^{\halfpos}-\mathcal{E}^{\halfneg}\big).
\end{equation}
In particular, if $(\vec f, \vec g)\equiv 0$, then the $\theta$-method conserves the discrete energy in the sense:
\begin{equation*}
\mathcal{E}^{\halfpos}=\mathcal{E}^{\halfneg}=\mathcal{E}^{\tfrac{1}{2}}.
\end{equation*}
\end{corollary}

\begin{proof}
We choose $(\vec \varphi_h, \vec \psi_h) = \urhdiff{n+1}{n-1}$ as a test function in~\eqref{eq:th} and obtain~\eqref{eq:nrg_estimate} with direct calculations.
\end{proof}

\subsubsection{Convergence}

To establish the convergence of the fully discrete approximation \cref{eq:th} to the elastic wave propagation problem \cref{eq:elasticwave}, we follow the arguments presented in \cite{Grote2009,Karaa}, resulting in the following convergence theorem.

\begin{theorem}[Error estimate]\label{thm:tau_h_bound} Assume that the solution pair $(\vec u, \vec r)$ to \cref{eq:elasticwave} satisfies the following regularity assumptions:
\begin{align*}
&(\vec u,\vec r) \in C^{2}([0,T]; \vec V^{p+1}\times \vec V^{p+1}),\\
&(\partial_t^3 \vec u, \partial_t^3\vec r) \in C([0,T]; L^2 \times L^2),\\
&(\partial_t^4 \vec u, \partial_t^4 \vec r) \in L^1(0,T; L^2 \times L^2),
\end{align*}
if $\theta \neq \tfrac{1}{12}$ and 
\begin{align*}
&(\vec u,\vec r) \in C^{2}([0,T]; \vec V^{p+1}\times \vec V^{p+1}),\\
&(\partial_t^5 \vec u, \partial_t^5\vec r) \in C([0,T]; L^2 \times L^2),\\
&(\partial_t^6 \vec u, \partial_t^6 \vec r) \in L^1(0,T; L^2 \times L^2),
\end{align*}
for $\theta = \tfrac{1}{12}$, where we have used the definition
\begin{align*}
      \vec V^k \coloneqq \left\{ \vec v \in \mathcal C^{0}(\Sigma)^3 \with \vec v|_\edge \in H^k(\edge)^3 ,\, \forall \edge \in \setEdge,\;\,  \vec v(\node) = \vec 0, \, \forall \node \in \setNodeDir \right\}.
 \end{align*}
Let the CFL condition \cref{eq:CFL} hold and assume the initial conditions are defined by~\cref{eq:intial_condition}.
Then, the fully discrete scheme given by \cref{eq:th}, where we use \cref{eq:central_diff_tau}, if $\theta \neq \frac{1}{12}$, and \cref{eq:central_diff_tau12} for $\theta = \frac{1}{12}$ to determine the fictitious values $(\vec u_h^{-1}, \vec r_h^{-1})$,
is uniquely defined in $\vec V_h \times \vec V_h$ and fulfills the $L^2$-error estimate
\begin{align*}
    \max_{0 \leq n \leq N} \| (\vec u_h^n, \vec r_h^n) - (\vec u(t_n), \vec r(t_n))\|_{L^2 \times L^2} \leq C_\mathrm{c} \big(h^{p+1} + \tau^s\big),
\end{align*}
with temporal order \( s = 2 \) if \( \theta \neq \tfrac{1}{12} \) and \( s = 4 \) if \( \theta = \tfrac{1}{12} \), where the constant $C_\mathrm{c}>0$ is independent of both \(h\) and \(\tau\), but depends linearly on \(T\).
\end{theorem}
\begin{proof}
The proof follows analogously to the arguments in \cite{Karaa}, where the $L^2$-scalar product on the left-hand side is replaced by the bilinear form $m$, the $L^2$-projection by the $m$-orthogonal projection and the discontinuous Galerkin space by the continuous finite element space defined in \cref{DGspace}.
\end{proof}

\begin{remark}[Alternative choices of projections]
    Assuming additional regularity on the initial data, the $m$-orthogonal projection $P_m$ can by replaced by a local projection or interpolant satisfying standard interpolation error bounds.
\end{remark}
 
\section{Subspace decomposition preconditioner}\label{sec:precond}

A practical implementation of the proposed method requires solving linear systems of equations arising from problem~\cref{eq:weakstationarydisc} in the stationary setting, or from problem~\cref{eq:th} (for  \(\theta>0\)) in the time-dependent setting. In engineering applications, these systems are typically large, sparse, and ill-conditioned. Efficient numerical solvers therefore rely on parallelization and on representations of the discretization across multiple scales. This principle underlies iterative methods such as geometric~\cite{Br77} and algebraic~\cite{XZ17} multigrid, as well as domain decomposition techniques~\cite{ToW05}. These approaches construct multilevel representations of the underlying graph (coarsening), which are then used to precondition the linear system. Spectral bounds for such multilevel preconditioners often yield proofs of optimal convergence rates. However, linear systems of equations originating from geometrically complex beam-network models typically require specially designed preconditioners that account for the intricacies of the geometry.

In the following, we consider a two-level overlapping additive Schwarz preconditioner, employing a direct solver for the local subproblems. This choice allows us to establish optimal convergence of the conjugate gradient method under mild assumptions on the underlying graph, following the analysis in~\cite{GoHeMa22}.
In this section, we focus on the case of linear finite element methods, i.e., $k = 1$ in \cref{sec:disc}. An integral part of the convergence proof of the preconditioned iteration is the formulation of a spectral equivalence result. To this end, we introduce two bilinear forms on the space $\vec{\mathcal{V}} \coloneqq \mathcal{V} \times \mathcal{V} \times \mathcal{V}$, where $\mathcal{V}$ denotes the space of real-valued functions on the node set $\setNode$ that vanish at Dirichlet nodes.
 Throughout, the value of a function $\vec \lambda \in \vec{\mathcal{V}}$ at a node $\node$ is denoted by $\vec \lambda_\node$. The first bilinear form corresponds to a mass-type operator, while the second defines a weighted graph Laplacian operator.
They are for any functions  $\vec \lambda, \vec \mu \in \vec{\mathcal{V}}$ defined~by
\begin{align}
	\label{EQ:graphlaplacian}
	\mathcal M(\vec \lambda,\vec \mu ) &\coloneqq \sum_{\node \in \setNode}\frac12\sum_{\edge \sim \node} \vec \lambda_\node \vec \mu_\node h_\edge,\\
	\mathcal L(\vec \lambda,\vec \mu) &\coloneqq \sum_{\node \in\setNode}\frac12 \sum_{\substack{\edge\sim \node\\ \edge =(\node,\node^\prime)}} \frac{(\vec \lambda_\node-\vec \lambda_{\node^\prime})\cdot (\vec \mu_\node-\vec \mu_{\node^\prime})}{h_\edge},
\end{align}
where the weighting with the edge length $h_\edge$ is chosen to be the same as for the mass and stiffness operators in a one-dimensional finite element implementation. 

\begin{lemma}[Spectral equivalence] \label{spectralequ}
Denote by \(\vec \xi = (\vec \xi_\node)_{\node \in \setNode}\) the vector of nodal values of the \(\mathcal{P}^1\) finite element function \(\vec u_h\), and similarly by \(\vec \eta = (\vec \eta_\node)_{\node \in \setNode}\) the vector of nodal values of the function \(\vec r_h\). Then, 
\begin{align}
\mathcal L(\vec \xi,\vec \xi) + \mathcal L(\vec \eta,\vec \eta) &\lesssim 	a\big((\vec u_h,\vec r_h),(\vec u_h,\vec r_h)\big) \lesssim \mathcal L(\vec \xi,\vec \xi) + \mathcal L(\vec \eta,\vec \eta),\label{eq:spectraleq1}\\
\mathcal M(\vec \xi,\vec \xi) + \mathcal M(\vec \eta,\vec \eta) &\lesssim m\big((\vec u_h,\vec r_h),(\vec u_h,\vec r_h)\big) \lesssim \mathcal M(\vec \xi,\vec \xi) + \mathcal M(\vec \eta,\vec \eta) .\label{eq:spectraleq2}
\end{align}
\end{lemma}

\begin{proof}
To prove estimate \cref{eq:spectraleq1}, we first note that the gradient of a linear finite element function is piecewise constant. Consequently, it follows that
\begin{align*} \|\vec u_h\|_{\vec V}^2 = \sum_{\edge \in \setEdge} \int_\edge |\partial_x \vec u_h|^2 \dx = \sum_{(\node,\node') = \edge \in \setEdge} \frac{(\vec \xi_{\node} - \vec \xi_{\node'})^2}{h_\edge} = \mathcal L(\vec \xi,\vec \xi). \end{align*}
Combining this identity with \cref{lem:coercivity} then yields the desired spectral equivalence between the stiffness operator and the component-wise graph Laplacian.

For the proof of estimate \cref{eq:spectraleq2}, we explicitly compute the relevant one-dimensional integrals, which gives
	\begin{align*}
		\|\vec u_h\|_{L^2}^2 = \sum_{\edge \in \setEdge} \int_\edge |\vec u_h|^2 \dx = \sum_{(\node,\node') = \edge \in \setEdge}\frac{h_\edge}{3}\big(\vec \xi_\node^2 + \vec \xi_\node \vec \xi_{\node'} + \vec \xi_{\node'}^2\big).
	\end{align*}
Applying Young's inequality, we obtain the following bounds for each summand:
	\begin{align*}
\frac{h_\edge}{6}\big(\vec \xi_\node^2 + \vec \xi_{\node'}^2\big) \leq \frac{h_\edge}{3}\big(\vec \xi_\node^2 + \vec \xi_\node \vec \xi_{\node'} + \vec \xi_{\node'}^2\big) \leq \frac{h_\edge}{2}\big(\vec \xi_\node^2 + \vec \xi_{\node'}^2\big).
	\end{align*}
As a direct consequence, we have
\begin{align*}
    \mathcal M(\vec \xi,\vec \xi) \lesssim \| \vec u_h \|_{L^2}^2 \lesssim \mathcal M(\vec \xi,\vec \xi),
\end{align*}
from which the desired spectral equivalence of the mass operator follows by \cref{eq:boundD}.
\end{proof}

\begin{remark}[Higher-order]
    The spectral equivalence extends to polynomial degrees beyond $k = 1$. For simplicity, however, we consider only the lowest-order~case.
\end{remark}

\begin{figure}[ht]
 \includegraphics[height=.375\textwidth]{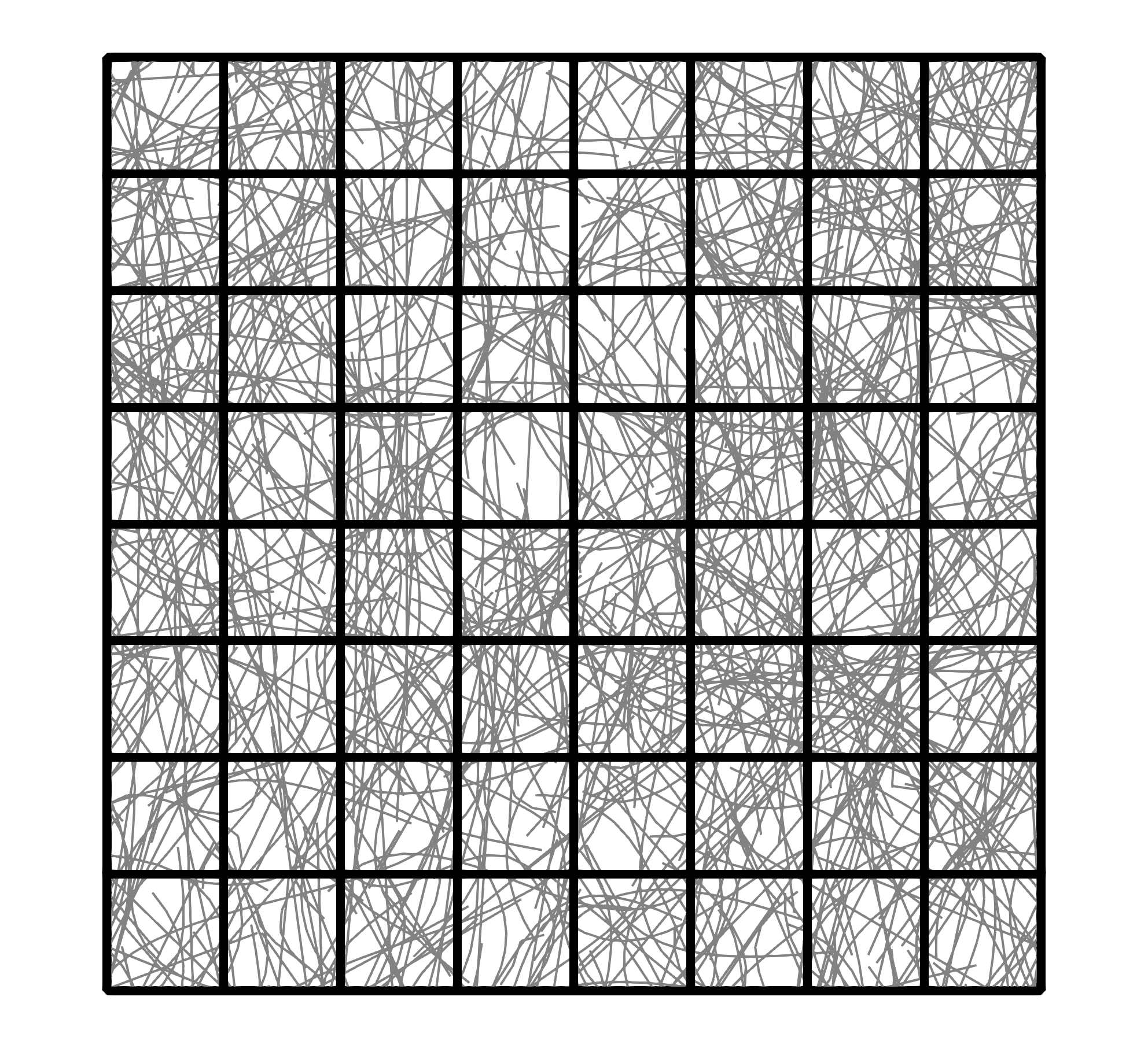}\hspace{.5cm}
 \includegraphics[height=.375\textwidth]{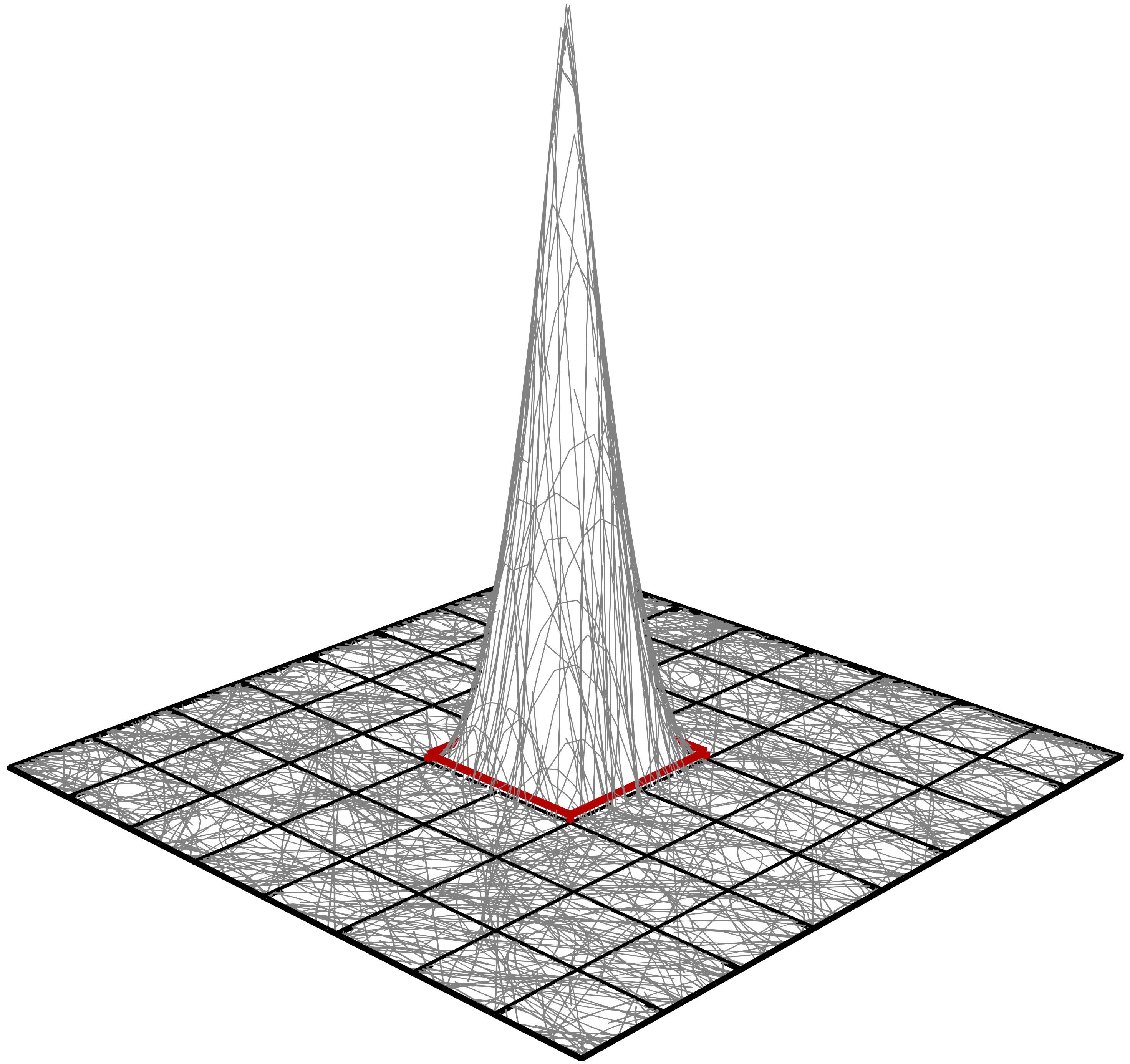}
 \caption{An artificial mesh $\mathcal{T}_H$ over a network (left) and a basis function~$\varphi_i$ with the boundary of its support marked in red (right), both in a two-dimensional setting.}\label{fig:artificialgrid}
\end{figure}

For the two-level overlapping additive Schwarz preconditioner, we introduce a coarse representation of the Timoshenko beam network on an artificial (coarse) mesh $\mathcal{T}_H$ over the bounded domain $\Omega \subset \mathbb{R}^3$, as illustrated in \cref{fig:artificialgrid} (left) for the two-dimensional case.
For simplicity, we assume that $\Omega$ is a box $[0,l_1] \times [0,l_2] \times [0,l_3]$ equipped with a uniform Cartesian mesh. Let $\{\varphi_i\}_{i=1}^m$ denote the set of trilinear basis functions, where $m$ is the number of nodes of the coarse mesh, with corresponding supports denoted by $U_i \coloneqq \operatorname{supp}(\varphi_i)$. We further define $\vec{U}_i \coloneqq U_i \times U_i \times U_i$. An illustration of a coarse-scale basis function in two dimensions is shown in \cref{fig:artificialgrid} (right). The space of continuous piecewise trilinear functions on $\mathcal{T}_H$ that satisfy Dirichlet boundary conditions on boundary segments where network nodes are fixed is denoted by $\mathcal{V}_H$. Note that the domain of a function in~$\mathcal{V}_H$ is considered to be the set of nodes of the spatial network.

We treat both the stationary problem and a single step of the time-dependent problem simultaneously by introducing a bilinear form $a^*$, which equals $a$ in the stationary case and $m + \tau^2 \theta a$ in the time-dependent case. To this end, we employ a preconditioner based on the following subspace decomposition:
\begin{equation*}
	\vec{\mathcal{V}}= \vec{\mathcal{V}}_{0}+   \vec{\mathcal{V}}_{1}+\dots+ \vec{\mathcal{V}}_{m}
\end{equation*}
with the coarse space $\vec{\mathcal{V}}_{0}\coloneqq \mathcal{V}_H\times \mathcal{V}_H \times \mathcal{V}_H$ and the local subspaces defined for $i=1,\dots,m$ by $\vec{\mathcal{V}}_{i}\coloneqq \{ \vec v\in \vec{\mathcal{V}} \with \operatorname{supp}(\vec v)\subset \vec U_i\}$. Given this  decomposition, we define, for each subspace, a corresponding subspace projection operator $P_i \colon \vec{\mathcal{V}} \times \vec{\mathcal{V}} \rightarrow \vec{\mathcal{V}}_i \times \vec{\mathcal{V}}_i$ satisfying,  for all $(\vec w, \vec z)\in \vec{\mathcal{V}}_{i}\times \vec{\mathcal{V}}_{i}$,
\begin{equation}
	\label{eq:subspcorr}
	a^*\big(P_i (\vec v, \vec s),(\vec w, \vec z)\big)=a^*\big((\vec v, \vec s),(\vec w, \vec z)\big).
\end{equation}
Note that the well-posedness of the subspace correction operator $P_i$ follows directly from the fact that $a^*$ defines an inner product on $\vec{\mathcal{V}} \times \vec{\mathcal{V}}$. We then define the preconditioned operator $P \colon \vec{\mathcal{V}} \times \vec{\mathcal{V}} \rightarrow \vec{\mathcal{V}} \times \vec{\mathcal{V}}$ as
\begin{equation*}
	P\coloneqq P_0+P_1+\dots+P_m.
\end{equation*}
The corresponding preconditioner, denoted by $B$, is defined via the relation $P = B A^*$, where $A^*$ is the matrix representation of the bilinear form $a^*$. We emphasize that the preconditioner $B$ is never explicitly formed; its application only requires evaluating the action of $P$. This, in turn, involves the direct solution of one coarse global problem and $m$ local problems of the form \cref{eq:subspcorr}. These subproblems are decoupled and can be solved in parallel. In the following, we employ this preconditioner within the conjugate gradient method.

We can establish optimal convergence of the preconditioned conjugate gradient method using the preconditioner $B$, provided that the coarse mesh size $H$ exceeds a critical length scale $R_0$. To formalize this, we introduce additional notation. Let us consider the boxes centered at $x = (x_1, x_2, x_3)$, defined as
\[
B_R(x) \coloneqq [x_1 - R,\, x_1 + R) \times [x_2 - R,\, x_2 + R) \times [x_3 - R,\, x_3 + R),
\]
where the half-open intervals are replaced by closed ones if \(x_i + R = l_i\). Moreover, let us denote the total length of all beams adjacent to the nodes in $B_R(x)$ by
\[
|1|^2_{\mathcal M, B_R(x)} \coloneqq 
\sum_{\node \in \setNode \cap B_R(x)} \frac{1}{2} \sum_{\edge \sim \node} h_\edge.
\]

We also denote by $\bar{d}_x$ the degree (i.e., the number of incident edges) of a node $x$ in a subgraph $\bar{\mathcal G} = (\bar{\mathcal N}, \bar{\mathcal E})$. For any subset of nodes $X \subset \bar{\mathcal N}$, we define its volume as
\[
\operatorname{vol}(X) \coloneqq \sum_{x \in X} \bar{d}_x,
\]
and we write $\bar{\mathcal E}(X, X') \subset \bar{\mathcal E}$ for the set of edges in the subgraph having one endpoint in $X$ and the other in $X'$.
The network assumptions under which we establish the convergence of the preconditioned conjugate gradient method can now be formulated as follows; see also \cite[Lem.~3.7]{GoHeMa22}.

 \begin{assumption}[Network assumptions]
    \label{ass:network}
    There is a length-scale $R_0$, a uniformity constant $\sigma$, and connectivity constants $\nu_1$ and $\nu_2$ 
     so that:
    \begin{enumerate}
    \item (homogeneity) for all $R \ge R_0$, 
      \begin{equation*}
 \max_{{B_R(x)}\subset \Omega} |1|^2_{\mathcal M,B_R(x)}\le \sigma \min_{{B_R(x)}\subset \Omega} |1|^2_{\mathcal M,B_R(x)};
      \end{equation*}
    \item (connectivity) for all $R \ge R_0$ and $x \in \Omega$, there is
      a connected subgraph of nodes
      $\mathcal{\bar{G}} = (\bar{\setNode}, \bar{\setEdge})$ (where $\bar{\setNode}$ denotes nodes and $\bar{\setEdge}$ edges connecting the nodes) of
      $\mathcal{G}$ containing
      \begin{enumerate}
      \item all edges with one or both endpoint in $B_R(x)$,
      \item only edges with endpoints contained in $B_{R+R_0}(x)$,
         \item 
      \begin{equation*}
        \text{vol}(\mathcal{\bar{N}}) \le \nu_1 \left(\frac{R}{R_0}\right)^d,
      \end{equation*}
      \item the following $d$-dimensional isoperimetric inequality holds,
      \begin{equation*}
        \left(\text{vol}(X)\right)^{(d-1)/d} \le \nu_2 |\mathcal{\bar{E}}(X, {\bar{\setNode}} \setminus X)|
      \end{equation*}
      for all $X \subset \mathcal{\bar{N}}$ assuming
      $\text{vol}(X) \le \text{vol}({\bar{\setNode}} \setminus X)$;
   \end{enumerate}
    \item (locality) for all edges in the network $\{x, y\} \in \mathcal{E}$, the edge length satisfies $|x - y| < R_0$;
    \item (boundary density) for every $y\in \setNode_D$, there is an 
      $x \in \setNode_D\setminus\{y\}$ such that
      $|x - y| < R_0$.
    \end{enumerate}
  \end{assumption}

The homogeneity assumption puts a bound on the density variation of the network on scales larger than $R_0$. The connectivity assumption requires that any subset of the network is well connected to the rest of the network. Under these two assumptions, the network behaves essentially as a continuous material on scales larger than $R_0$. The locality condition can be ensured by further discretization, and the boundary density assumption guarantees the presence of a sufficient number of fixed nodes on the constrained boundary segments.

We now apply \cite[Thm.~4.4]{GoHeMa22} to show convergence of the preconditioned conjugate gradient iteration. This result requires a sparse, symmetric system matrix and that the network satisfies \cref{ass:network}. Additionally, the system matrix must be spectrally equivalent to the reciprocal edge-weighted graph Laplacian $\mathcal L$. In the stationary case, this follows directly from \cref{spectralequ}. For the wave equation, we follow the argument of \cite[Lem.~4.2]{GoHeMa22} with $\mathcal{L}$ replaced by $\mathcal{M} + \tau^2 \theta \mathcal{L}$, again using the spectral equivalence from \cref{spectralequ}. In both cases, we obtain the same optimal convergence result, with condition number $\kappa$ independent of $H$ and $h$.

\begin{theorem}[Preconditioned conjugate gradient method]\label{thm:cg}
If $H \geq 2 R_0$, the preconditioned conjugate gradient iterate $(\vec u_h^{(\ell)}, \vec r_h^{(\ell)})$ satisfies
\[
\|(\vec u_h, \vec r_h) - (\vec u_h^{(\ell)}, \vec r_h^{(\ell)})\|_{\vec V \times \vec V} 
\leq 2 \left( \frac{\sqrt{\kappa}-1}{\sqrt{\kappa}+1} \right)^\ell 
\|(\vec u_h, \vec r_h) - (\vec u_h^{(0)}, \vec r_h^{(0)})\|_{\vec V \times \vec V},
\]
where the condition number $\kappa$ of the preconditioned operator depends only on the constants in \cref{spectralequ} and on $\nu_1$ and $\nu_2$ from the network assumptions.
  \end{theorem}

The proof of the theorem is inspired by classical Schwarz theory (see, e.g.,~\cite{Xu92,KorY16}) and involves constructing a quasi-interpolation operator in the spatial network setting, whose approximation and stability properties can be established using Poincar\'e and Friedrichs inequalities on subgraphs. In practice, this preconditioner has proven effective in handling the typically complex geometry of spatial networks and highly varying material properties; see, e.g.,~\cite{Grtz2024}.

\begin{remark}[Time-dependent case]
Note that in the time-dependent case, all occurrences of $(\vec u_h, \vec r_h)$ in the preceding section should be understood as $(\vec u_h^{\,n}, \vec r_h^{\,n})$. The solution from the previous time step provides a natural initial guess for the iteration. Moreover, if the time step and the material parameters are constant in time, the same matrices appear in the local problems at each step, allowing their factorizations to be saved and reused to further accelerate the computations.
\end{remark}

\section{Numerical examples}\label{sec:numexp}

We consider two engineering applications to verify our theoretical findings.

\subsection{Elliptic $h$-convergence}
\begin{figure}
\centering
\includegraphics[width = 0.60
\textwidth]{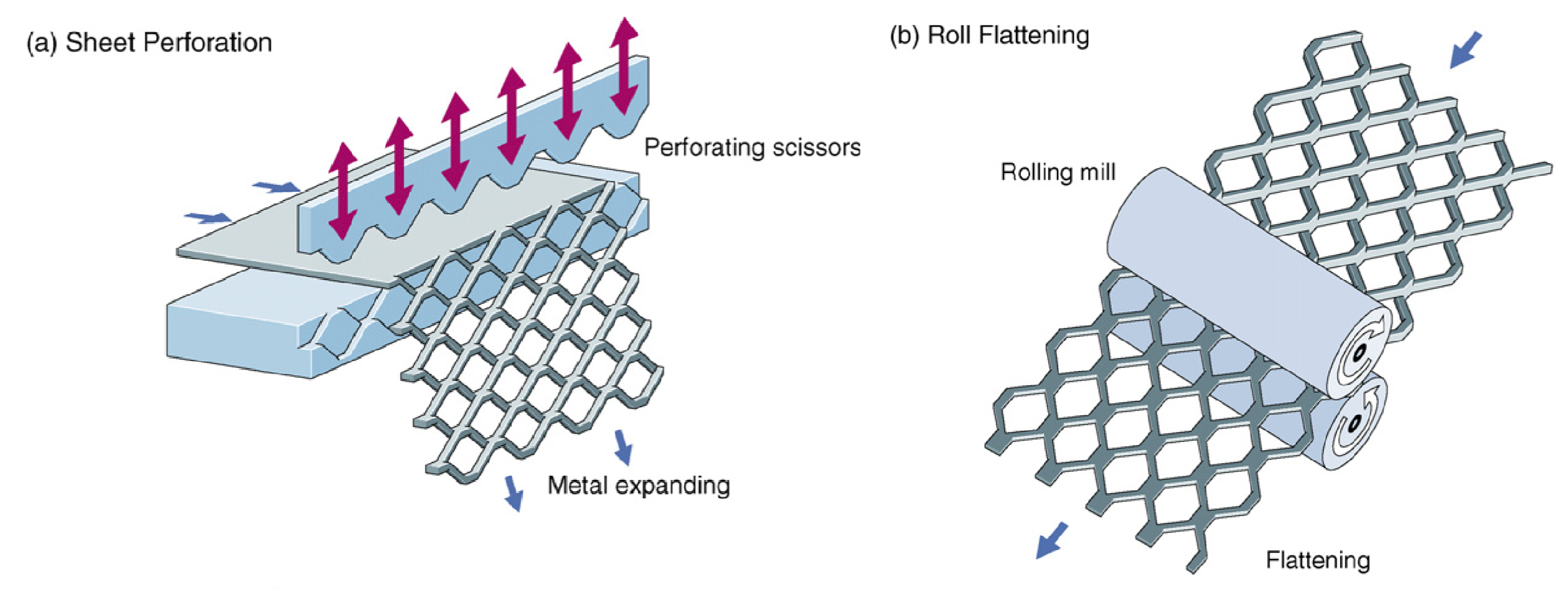}
\hspace{0.05\textwidth}
\includegraphics[width = 0.25\textwidth]{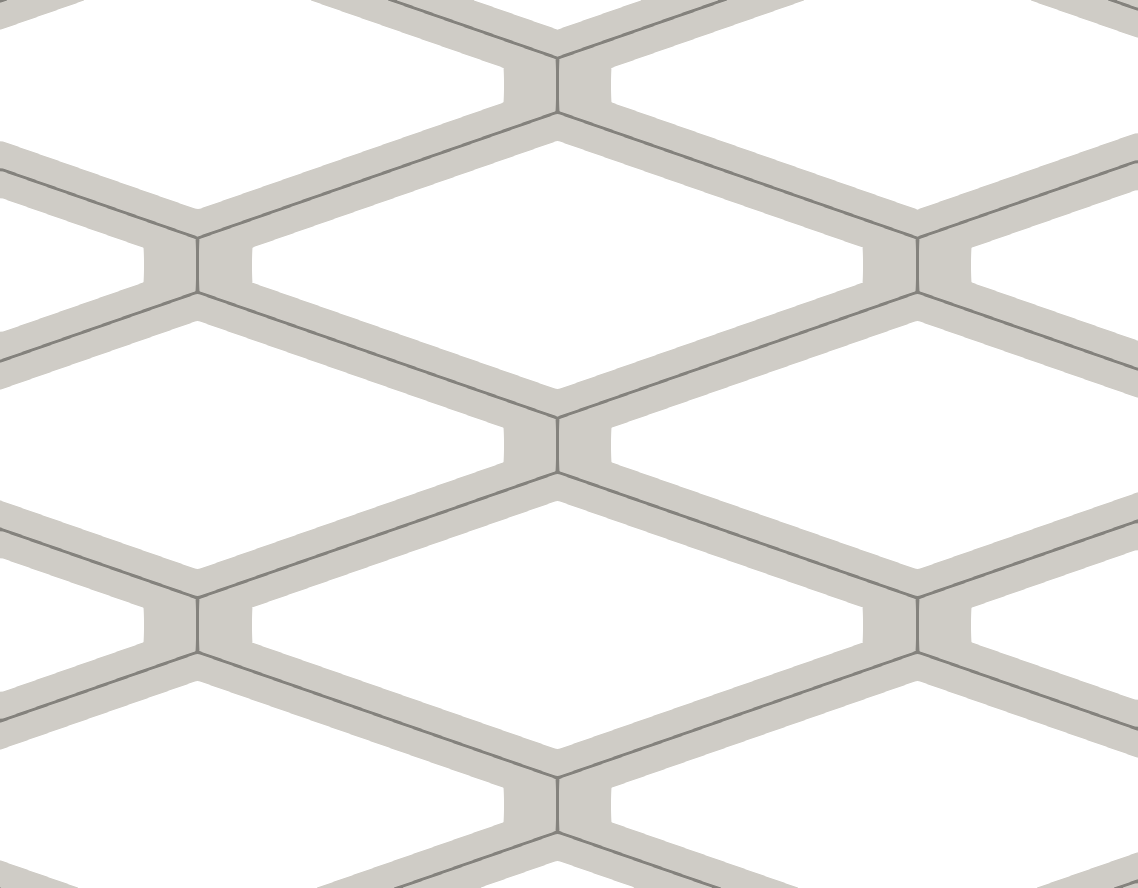}
\caption{The first two figures from the left are illustrations of the production of expanded metal, as well as the optional roll flattening post processes   \cite{ExpandedMetalIllustration}. The right figure presents the network structure of a two-dimensional rolled expanded metal. }
\label{fig:expanded}
\end{figure}

In the first numerical example, we investigate the theoretical error bound from \cref{thm:conv_stat} by simulating a small piece of expanded metal. Expanded metal is a mesh material produced by perforating a metal sheet. \cref{fig:expanded} shows illustrations of the manufacturing process and the corresponding network model. The large diamonds in the model measure 80~mm in width and 40~mm in height, while the slender mesh strands are 6~mm wide and 3~mm thick, and are modeled as cold-rolled steel (isotropic, $E = 210$~GPa, $\nu = 0.3$, $\rho = 7850$~kg/m$^3$).

To analyze the elliptic error bound, we consider a $32$~cm $\times$ $32$~cm metal mesh consisting of $4 \times 8$ diamonds. The network model is generated in the $xy$-plane, with the coarsest network, $\mathcal{G}_0$, containing the minimal number of nodes required to capture the geometry. The longest edge in this network, $h_0$, is approximately $45$~mm. Higher-resolution models are obtained by subdividing the edges of the coarse network. Specifically, we consider networks $\mathcal{G}_i$ with minimal components such that the largest edge length is $h_i = h_0 / 2^i$, for $i = 0, 1, \dots, 6$, and $i = 10$.
The mesh is bent out of the plane by clamping the right boundary and imposing Dirichlet conditions only on the displacement in the $z$-direction on the opposite side to achieve roughly a $15^\circ$ bend $(\text{tan}(15^\circ) \,\cdot\, $32 cm). The elliptic solutions $(\vec u_{h_i}, \vec r_{h_i})$ are then found by solving the corresponding linear systems of equations for each of the networks $\mathcal G_i$. The coefficients in the linear system of equations are based on the units (metric ton, cm, s).
The solution of the detailed model, $(\vec u_h, \vec r_h) \coloneqq (\vec u^{10}_{h_{10}}, \vec r^{10}_{h_{10}})$, is treated as the reference solution. The coarser solutions are upscaled via linear interpolation, consistent with the $\mathcal{P}^1$ finite element basis used in the discretization, and then compared to the reference.

\begin{figure}
\centering
\includegraphics[width = 0.4
\textwidth]{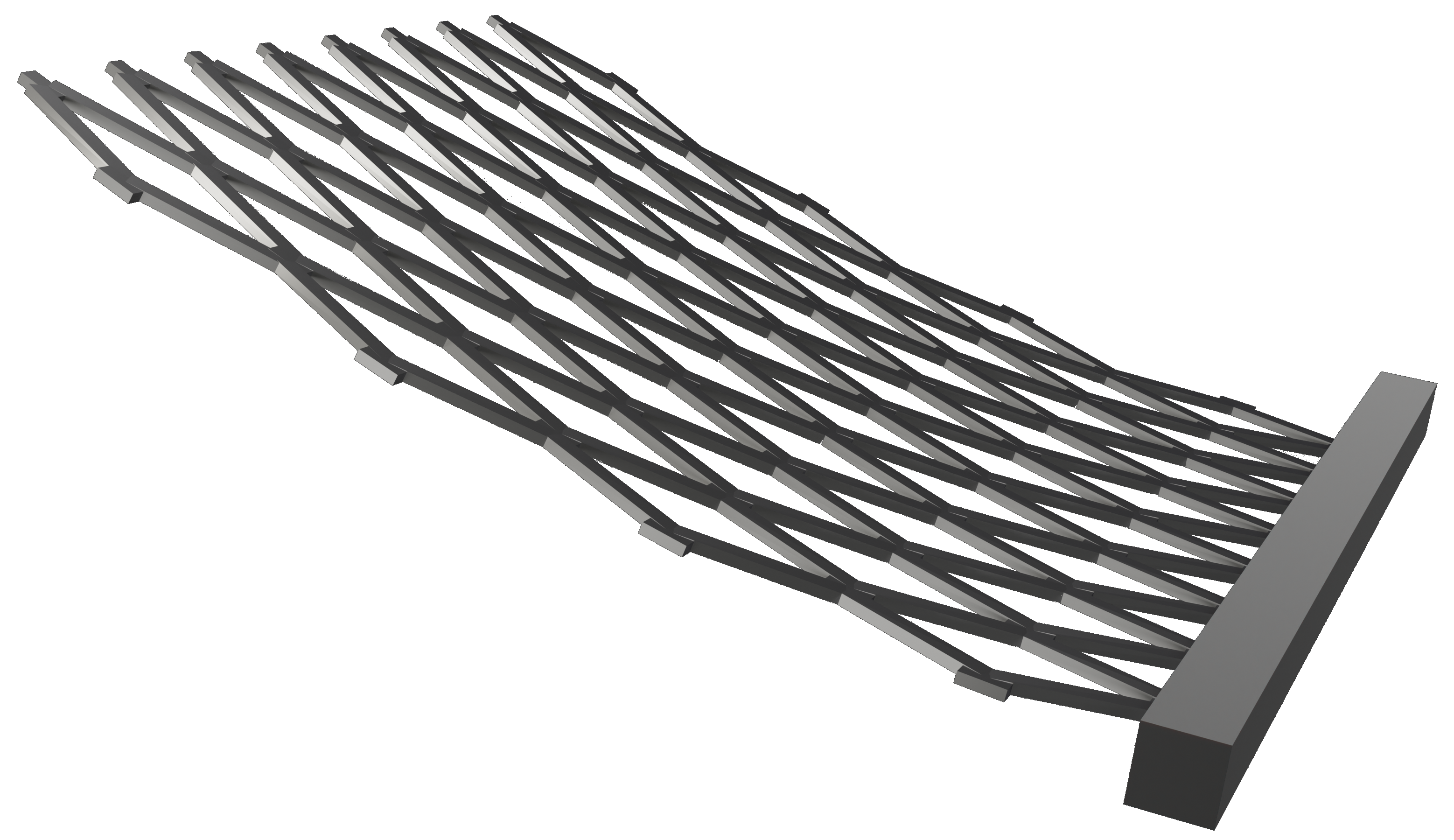}
\includegraphics[width = 0.35\textwidth]{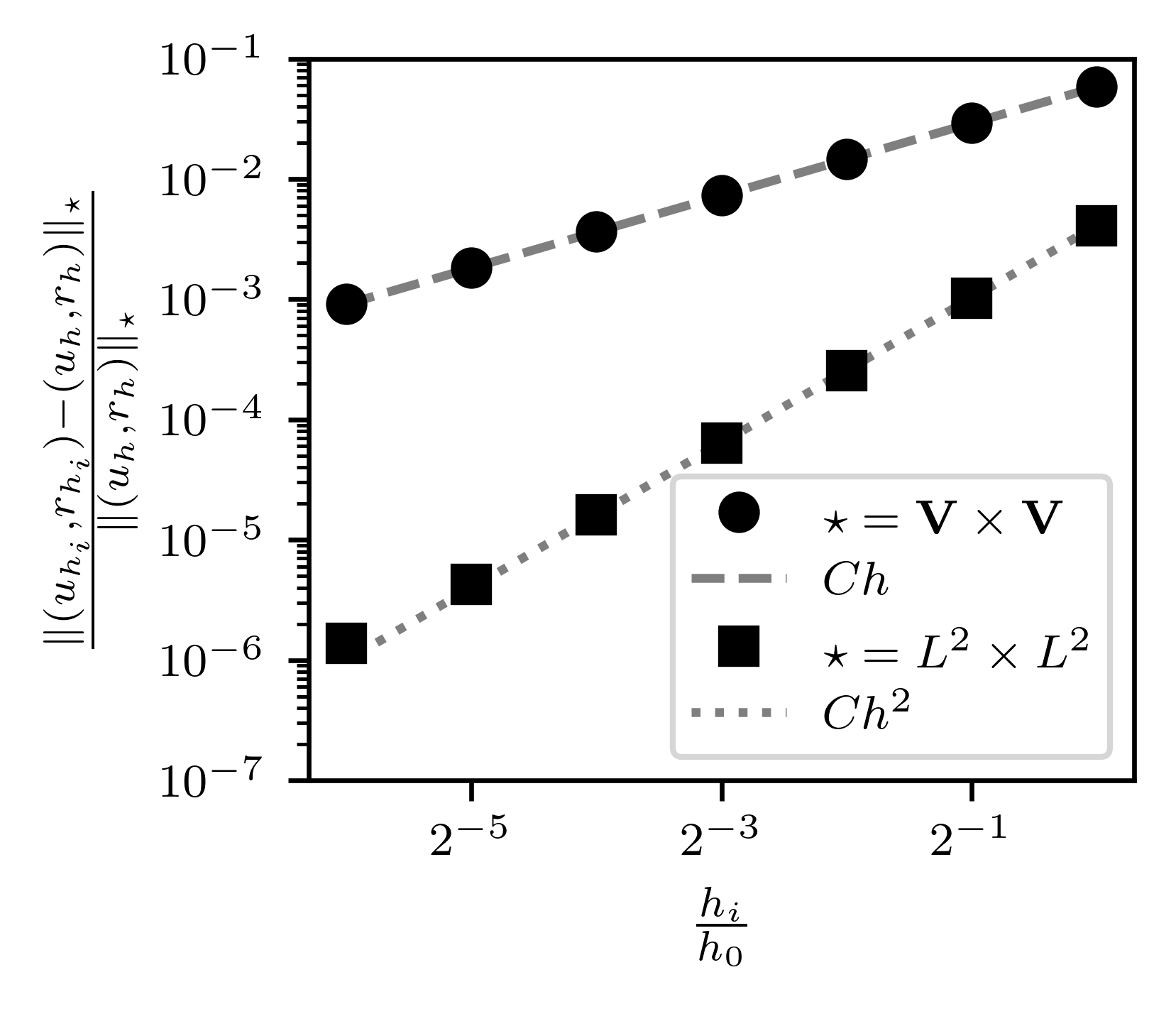}
\caption{The left figure shows the elliptic solution for $\mathcal{G}_{10}$, while the right figure displays the  convergence of the method with respect to the spatial discretization parameter $h_i$.
}
\label{fig:numeric_spatial_elliptic}
\end{figure}

Illustrations of the elliptic solution and the error convergence are shown in \cref{fig:numeric_spatial_elliptic}. The results clearly exhibit the theoretical $h$-scaling in the $\mathbf{V} \times \mathbf{V}$ norm predicted by \cref{thm:conv_stat}, as well as $h^2$-scaling with respect to the $L^2 \times L^2$ norm.

\subsection{Numerical estimate of the CFL condition}
\begin{figure}
\centering
\includegraphics[width = 0.4\textwidth]{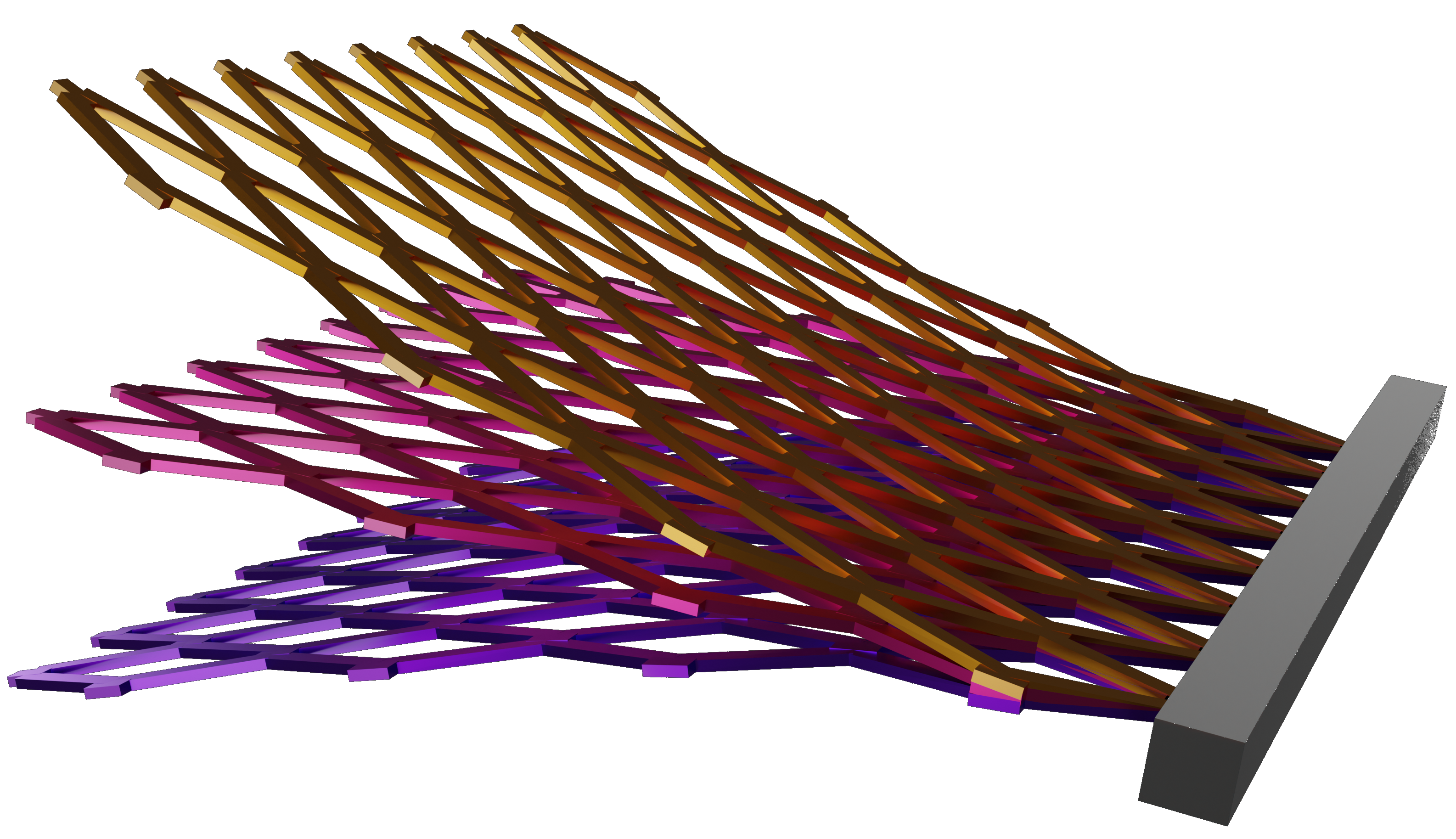}
\hspace{0.05\textwidth}
\includegraphics[width = 0.4\textwidth]{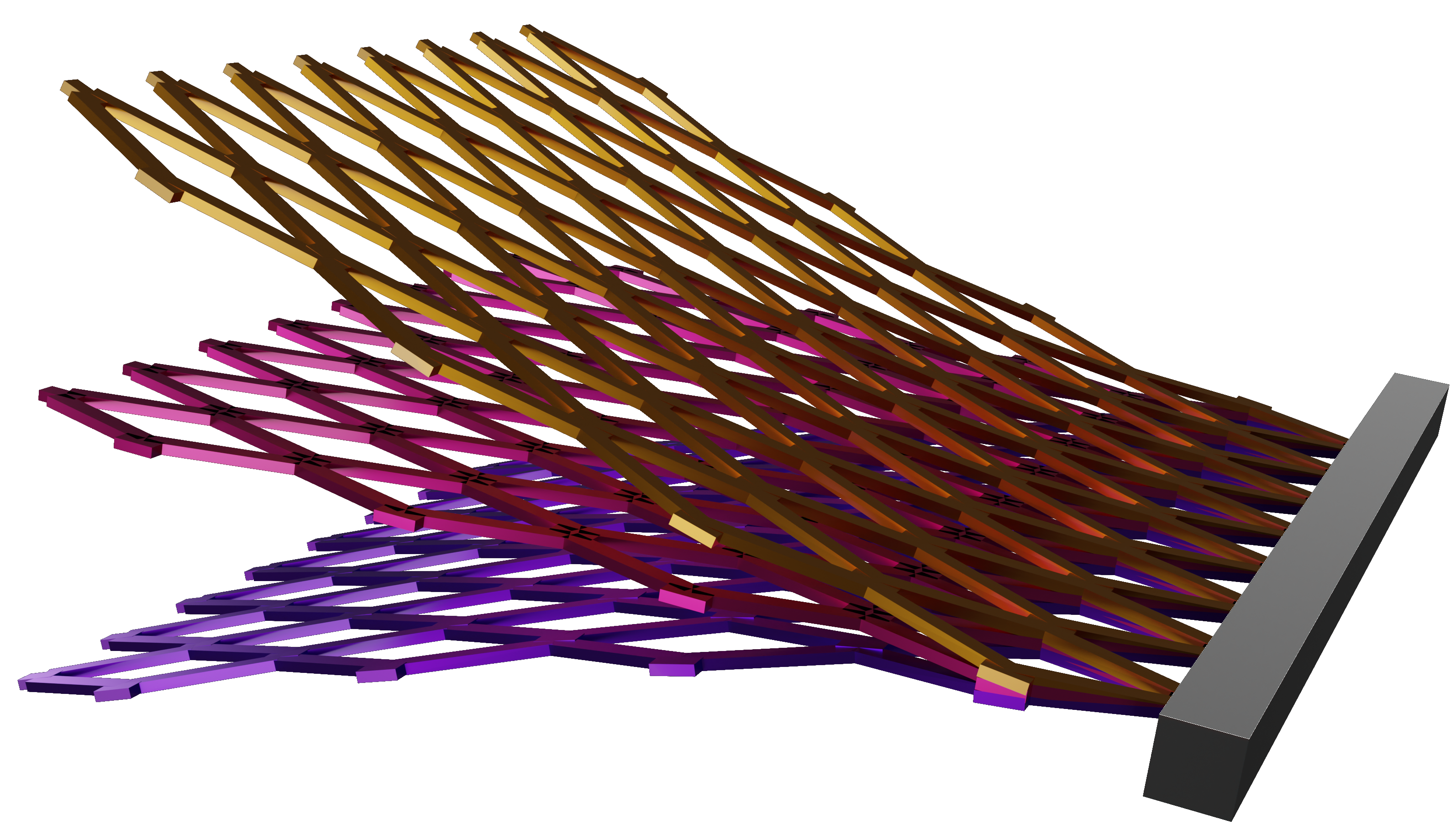}
\caption{The solutions of the two (similar) elastic waves are shown at approximately $t = 0$~ms (yellow), $7$~ms (pink), and $14$~ms (purple). The left figure is  the response obtained by removing the displacement condition in \cref{fig:numeric_spatial_elliptic}, while the right figure shows the first (scaled) eigenmode with the clamped boundary condition.
}
\label{fig:expanded-wave}
\end{figure}

In the second numerical example, we study elastic wave propagation in a network, focusing on the CFL condition for the leapfrog method ($\theta = 0$). The condition is evaluated using the coarse network model $\mathcal{G}_0$ from the previous example. The initial values $(\vec u_{h}^0, \vec r_{h}^0)$ are taken from the elliptic solution in \cref{fig:numeric_spatial_elliptic}, and we set $(\vec v_{0}, \vec p_{0}) = (0,0)$ in \eqref{eq:central_diff_tau}.

The specific wave problem considered is the response when the displacement conditions are released. This response is computed using the leapfrog method, starting with $(\vec u^0_{h}, \vec r^0_{h}) = (\vec u_{h_0}, \vec r_{h_0})$, while $(\vec u^1_h, \vec r^1_h)$ is obtained using \eqref{eq:th}, solving
$$  \frac{2}{\tau^2}m\big( ( \vec u_h^{1}, \vec r_h^{1}), (\vec \varphi_h, \vec \psi_h)\big)  
 = \frac{2}{\tau^2}m\big( ( \vec u_h^{0}, \vec r_h^{0}), (\vec \varphi_h, \vec \psi_h)\big) -a\big( (\vec u^0_h, \vec r^0_h), (\vec \varphi_h, \vec \psi_h) \big) $$
for all $(\vec \varphi_h, \vec \psi_h) \in \vec V_h \times \vec V_h$. \cref{fig:expanded-wave} shows the state of the mesh at various times after the displacement condition is released. In this example, roughly one periodic cycle is simulated ($T = 28$~ms).

The CFL condition for the leapfrog method in \eqref{eq:CFL} is evaluated as
\[
C_\text{CFL} = \sqrt{\frac{a_\text{max}}{m_\text{min}} \cdot \frac{1}{4}}   C_\text{inv}, \qquad
\tau< \frac{h_{\text{min}}}{C_\text{CFL}} = 41 \text{ ns},
\]
where $C_\text{inv} = \sqrt{12}$ in the present setting of linear finite elements in one dimension, and $h_{\text{min}}$ denotes the length of the smallest edge. Based on this bound, the leapfrog method was tested for $\tau$ near the CFL condition, and the $\tau$ at which the method became unstable was isolated. The results are presented in \cref{tab:leapfrog} and show that the instability occurs for $\tau$ roughly 7 times larger than the theoretical bound in \cref{thm:stability}. The unconditionally stable Crank–Nicolson method produced stable results for all of these cases, as well as for much larger $\tau$, confirming \cref{thm:stability}.

\begin{table}
\centering
\caption{Stability analysis of the leapfrog method for different time steps.
}
\begin{tabular}{c || c | c | c } 
 $\tau$ &  287 ns & 297 ns & 308 ns\\
 Time steps & 98 000  & 94 000 & 91 000 \\
 Stable & Yes  & Yes & No \\
\end{tabular}

\label{tab:leapfrog}
\end{table}

\subsection{Crank--Nicolson $\tau$-convergence}
Next, the time discretization error is assessed by comparing iterations for various $\tau$ with the Crank--Nicolson method, using a known semi-discrete solution. The network model is the same as in the previous two examples, and a semi-discrete problem is constructed using the first eigenmode of the corresponding elliptic operator. Specifically, we seek $(\hat{\vec{\phi}}_h^1, \hat{\vec{\nu}}_h^1) \in \vec V_h \times \vec V_h$ associated with the smallest (positive) eigenvalue $\lambda$ such that
\[
a\big((\hat{\vec{\phi}}_h^1, \hat{\vec{\nu}}_h^1), (\vec \varphi_h, \vec \psi_h)\big) = \lambda \, m\big((\hat{\vec{\phi}}_h^1, \hat{\vec{\nu}}_h^1), (\vec \varphi_h, \vec \psi_h)\big), \quad \forall (\vec \varphi_h, \vec \psi_h) \in \vec V_h \times \vec V_h.
\]
The eigenvector $(\hat{\vec{\phi}}_h^1, \hat{\vec{\nu}}_h^1)$ is computed numerically and, after rescaling, closely resembles the steady-state solution in \cref{fig:numeric_spatial_elliptic}. The smallest eigenvalue, approximately $\lambda_1 \approx 0.0584$, determines the resonance period $2\pi / \sqrt{0.0584} \ \text{ms} \approx 26 \ \text{ms}$, which agrees well with the period observed above in the CFL example (28 ms).

The semi-discrete solution we consider is
\[
(\vec u_h(t), \vec r_h(t)) = (\hat{\vec{\phi}}_h^1, \hat{\vec{\nu}}_h^1) \cos(\sqrt{\lambda_1}\, t).
\]
This function satisfies the semi-discrete problem
\begin{subequations}
\label{eq:wave-tau-conv}
\begin{equation}
m\big( (\ddot{\vec u}_h, \ddot{\vec r}_h), (\vec \varphi_h, \vec \psi_h)\big)
+ a\big( (\vec u_h, \vec r_h), (\vec \varphi_h, \vec \psi_h)\big)
= 0,
\quad \forall (\vec \varphi_h, \vec \psi_h) \in \vec V_h \times \vec V_h,
\end{equation}
with the initial conditions
\begin{equation}
(\vec u_h(0), \vec r_h(0)) = (\hat{\vec{\phi}}_h^1, \hat{\vec{\nu}}_h^1), 
\qquad
(\dot{\vec u}_h(0), \dot{\vec r}_h(0)) = (0,0).
\end{equation}
\end{subequations}

This problem is then solved using the Crank--Nicolson method, with $(\vec u_h^0, \vec r_h^0) = (\vec u_h(0), \vec r_h(0))$, and the first time step $(\vec u_h^1, \vec r_h^1)$ computed by solving
\begin{align*}
\frac{2}{\tau^2}
&m\big( (\vec u_h^{1}, \vec r_h^{1}), (\vec \varphi_h, \vec \psi_h)\big)
+ 
\frac{1}{2}a\big( (\vec u_h^{1}, \vec r_h^{1}), (\vec \varphi_h, \vec \psi_h)\big)
\\[2pt]
&=
\frac{2}{\tau^2}
m\big( (\vec u_h^{0}, \vec r_h^{0}), (\vec \varphi_h, \vec \psi_h)\big)
-
\frac{1}{2}a\big( (\vec u_h^{0}, \vec r_h^{0}), (\vec \varphi_h, \vec \psi_h)\big)
\end{align*}
for all $(\vec \varphi_h, \vec \psi_h) \in \vec V_h \times \vec V_h$. The latter problem corresponds to \eqref{eq:th} with zero initial velocities 
$(\vec v_{0}, \vec p_{0}) = (0,0)$ in \eqref{eq:central_diff_tau}.

\begin{figure}
\centering
\includegraphics[width = 0.35\textwidth]{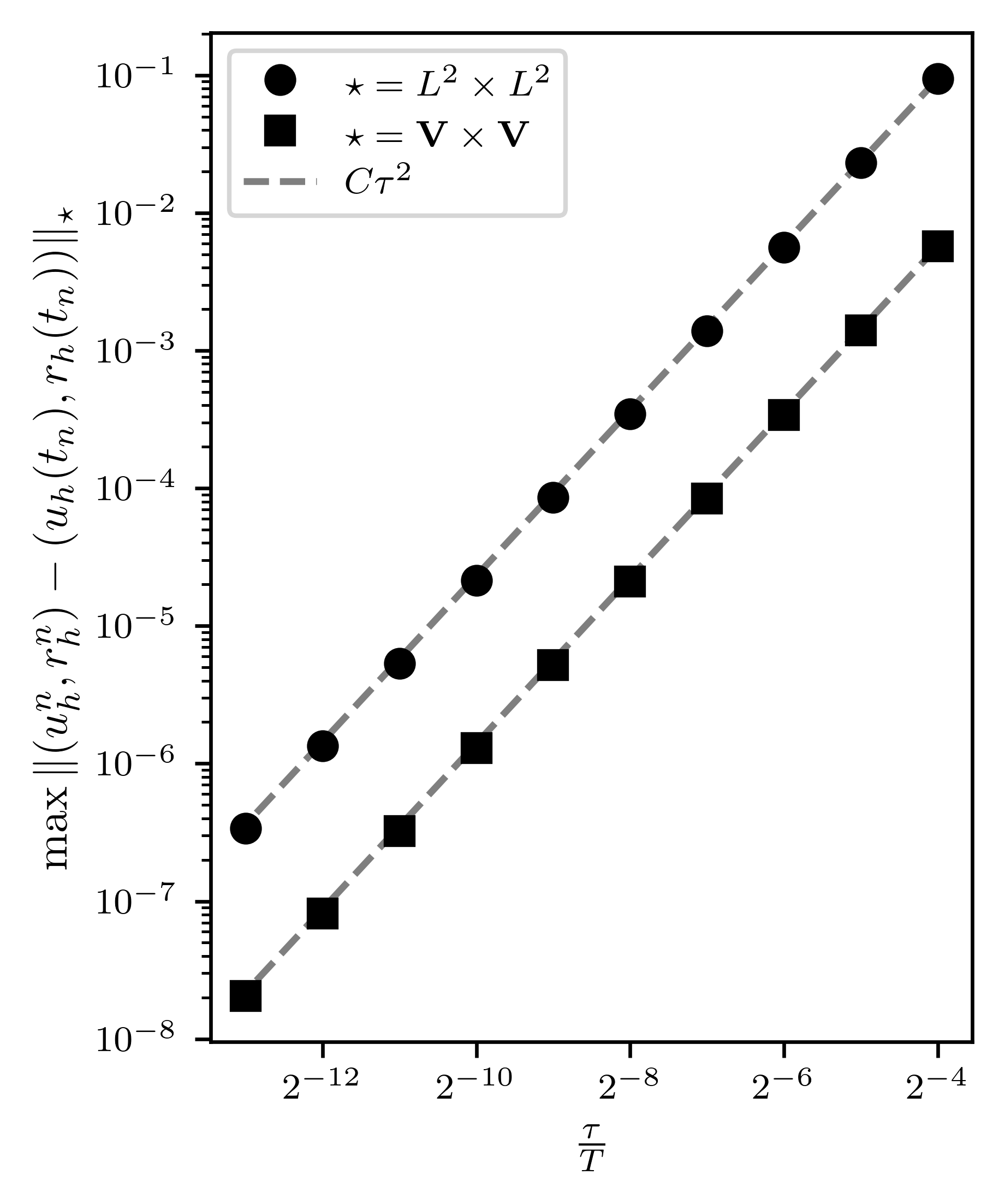}\hspace{.25cm}
\includegraphics[width = 0.35\textwidth]{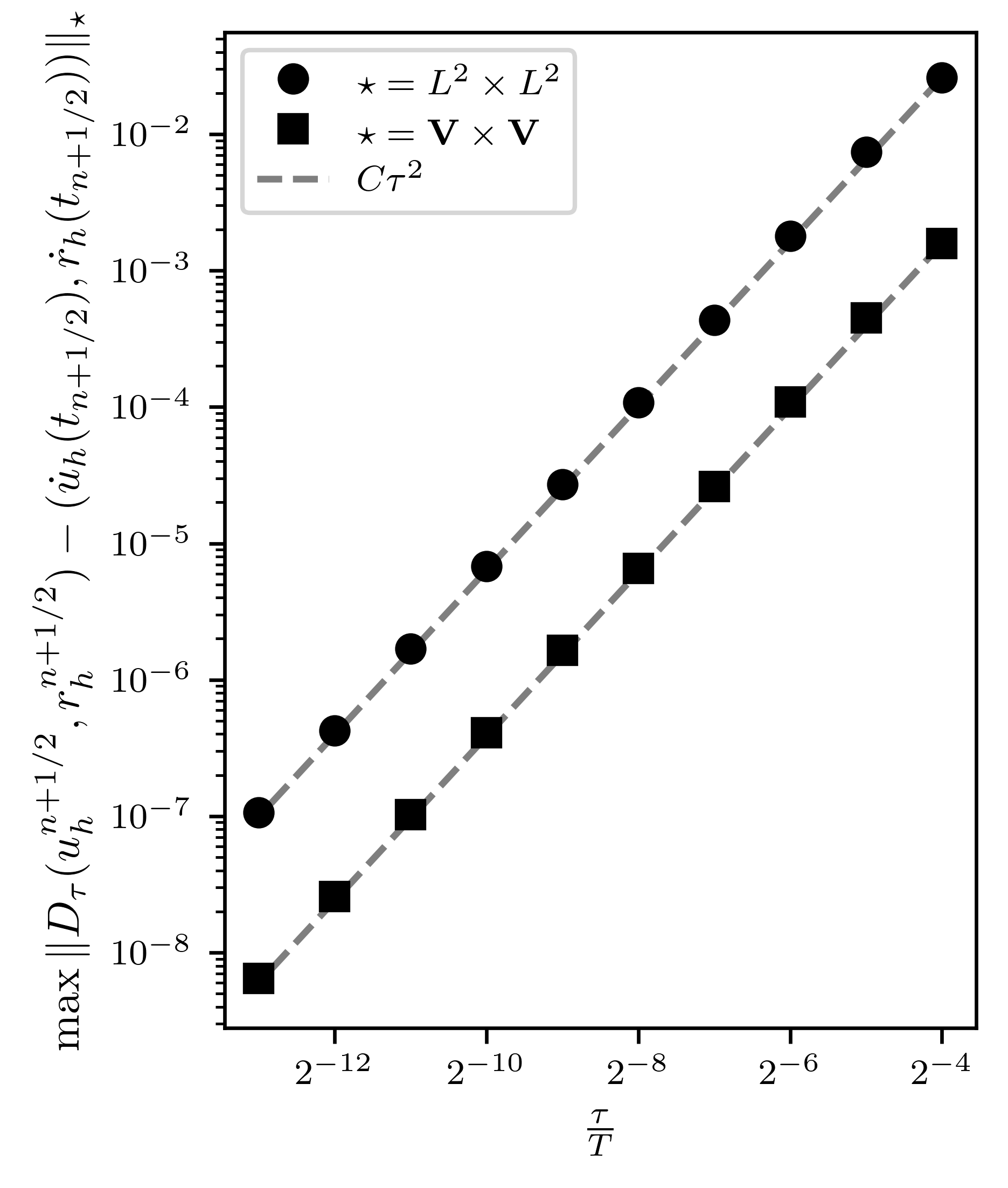}
\caption{Temporal convergence of Crank--Nicolson for a wave problem with a semi-discrete analytical solution.
}
\label{fig:tau-conv}
\end{figure}

We evaluate the temporal discretization error for various $\tau$ via
\begin{align*}
    &\max_{1 \leq n \leq N} \| (\vec u_h^n, \vec r_h^n) - (\vec u_h(t_n), \vec r_h(t_n)) \|,\\
&\max_{1 \leq n \leq N} \|  (D_\tau\vec u_h^{\,n+1/2}, D_\tau\vec r_h^{\,n+1/2}) - (\dot{\vec u}_h(t_{n+1/2}), \dot{\vec r}_h(t_{n+1/2})) \|,
\end{align*}
both measured in the $\mathbf{V} \times \mathbf{V}$ norm and the $L^2 \times L^2$ norm. The results, which numerically confirm the temporal convergence rates of the Crank--Nicolson method predicted in \cref{thm:tau_h_bound}, are shown in \cref{fig:tau-conv}.

\subsection{Paper example}
In the fourth numerical example, we consider a small piece of printing paper with a typical surface weight of $70~\text{g}/\text{m}^2$. The material is modeled by generating paper fibers (unbleached kraft) from geometrical distributions obtained via image analysis, including fiber lengths, widths, and curvature \cite{Grtz2024}. In \cite{Grtz2024}, the mechanical properties of paperboard (tensile stiffness, bending resistance) were studied using a network model with Timoshenko beams, with some simulations requiring the domain decomposition method from \cite{GoHeMa22}. The Timoshenko beam model is appropriate here because the fibers are densely packed, making the beam length (between intersections) small relative to the total fiber length.
Paper is an orthotropic material, being substantially stiffer in one direction due to the papermaking process. This direction is called the Machine Direction (MD), while the perpendicular Cross Direction (CD) is least stiff. In this example, the paper model captures orthotropy through biased fiber orientations. The following simulations consider only the machine direction; details on fiber properties and cross-direction simulations are provided in \cite{Grtz2024}.

\begin{figure}
\centering
\includegraphics[width = 0.95\textwidth]{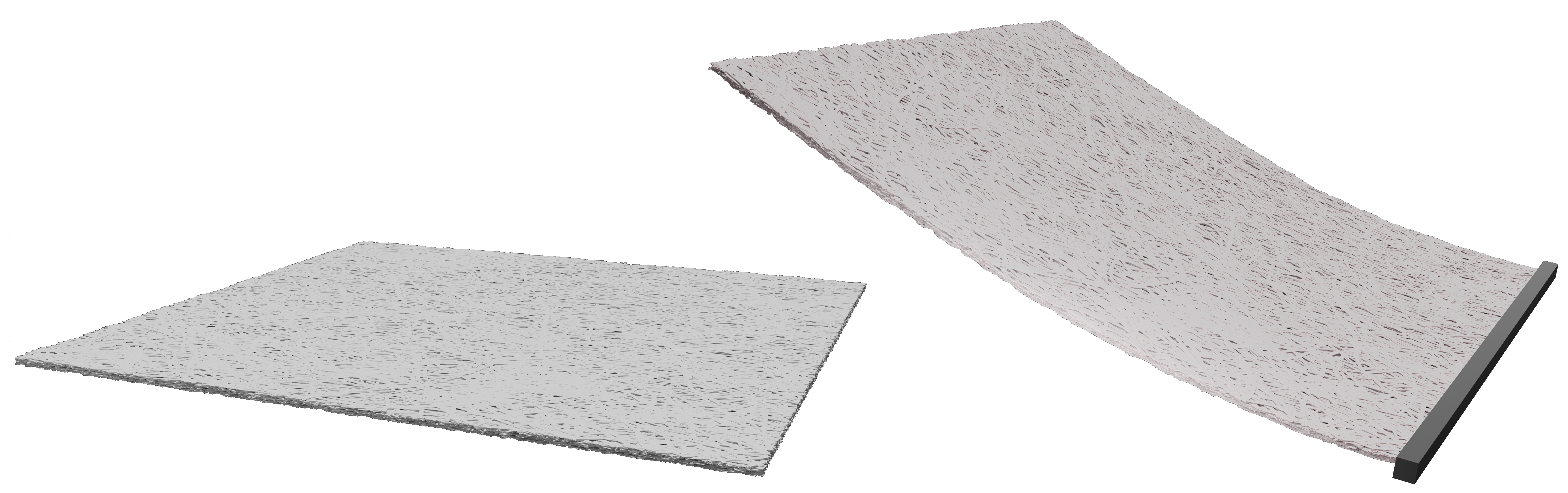}
\caption{The left figure shows the $8~\text{mm} \times 8~\text{mm}$ paper model, while the right figure displays the solution of an elliptic problem simulating a two-point bending experiment. }
\label{fig:fiber_model}
\end{figure}

The problem considered models an $8~\text{mm} \times 8~\text{mm}$ piece of paper, composed of 13,630 individual fibers. These fibers are initially discretized into $0.1~\text{mm}$ segments, with cross-sections derived from experimental data. The fiber cross-sections are either hollow (elliptic, 3~\textmu m wall thickness) or collapsed (rectangular). After deposition, the network consists of thousands of disconnected chains of edges. To form a connected network, all contacts are identified based on the 3D volume of the fibers. Beams are placed between each pair of contacting fibers by further discretizing the fibers at the closest points. The beams have rectangular cross-sections, determined by projecting the intersection outline onto a contact plane and fitting an appropriately oriented rectangle. Finally, any remaining disconnected edges are removed. The resulting network is well connected with small density variation on the scales $H$ used in the subspace decomposition for the numerical simulations. A visualization of the considered fiber network model is shown in \cref{fig:fiber_model}, and some general properties are summarized in \cref{tab:model-paper-prop}.

\begin{table}
\centering
\caption{Model properties}
\label{tab:model-paper-prop}
\begin{tabular}{l|r}  
Surface weight & 70 g/$\text{m}^2$  \\
Density ($\rho_s$)   & 770 kg/$\text{m}^3$ \\
Thickness ($t_s$) [ Surface weight/$\rho_s$] & 91 \textmu m \\
Model tensile stiff. (MD,CD) &  (12.0, 2.80)~GPa, \\
Num. beams & 1.45 million\\
Num. nodes & 0.808 million \\
Degrees of freedom & 4.85 million
\end{tabular}
\end{table}

\subsubsection*{Optimality of the preconditioner with respect to $H$}
An elliptic problem is considered to evaluate the $H$-invariance of the domain decomposition method in terms of its number of iterations. The setup is similar to \cref{fig:numeric_spatial_elliptic}, with one side clamped and the other bent out of plane by $15^\circ$ through $z$-directional Dirichlet conditions $(\tan(15^\circ) \cdot 8~\text{mm})$. The steady-state solution is computed in units (kg, \textmu m, ms) using a direct linear solver and shown in \cref{fig:fiber_model}.

\begin{figure}
\centering
\includegraphics{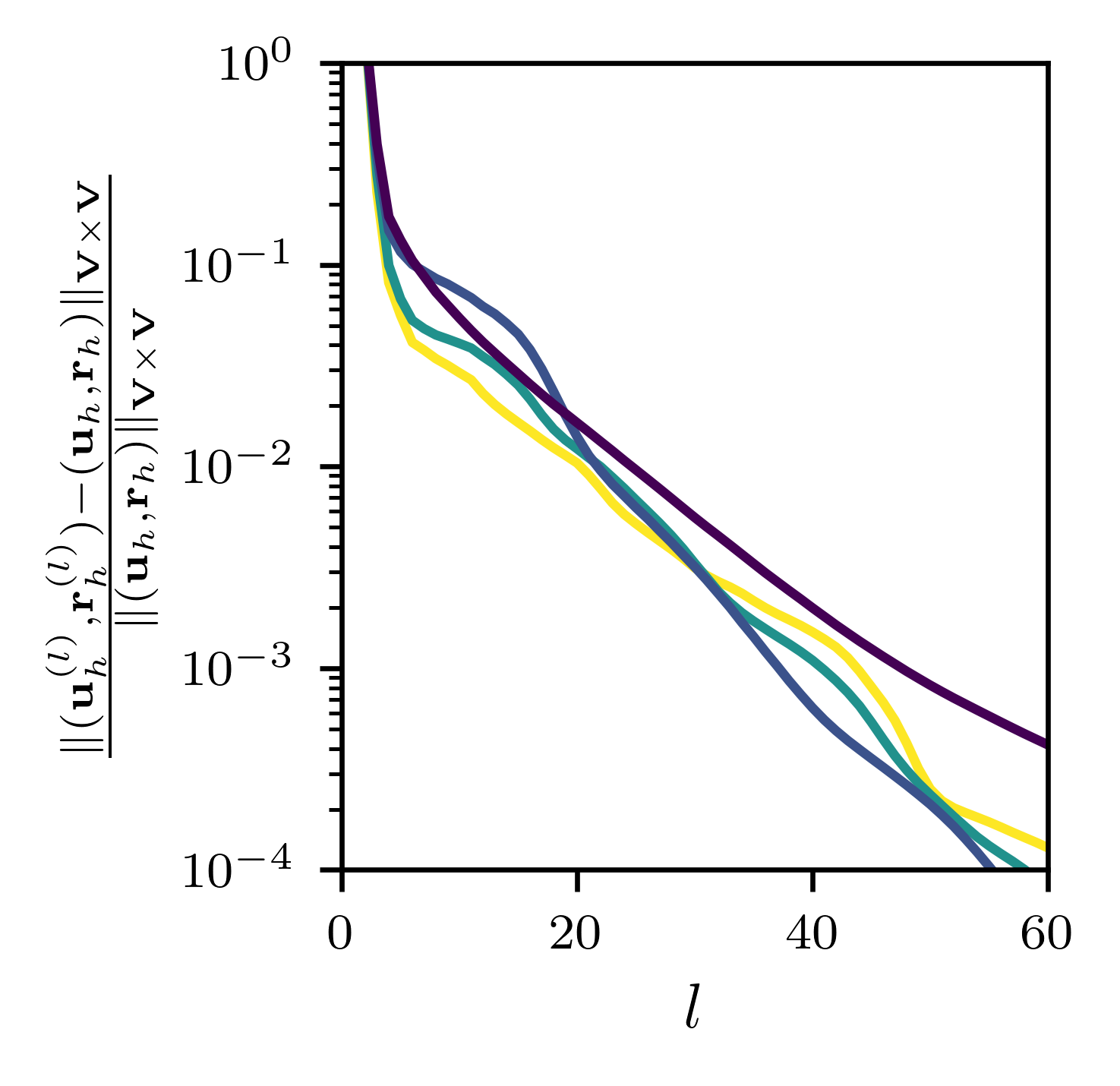}
\includegraphics{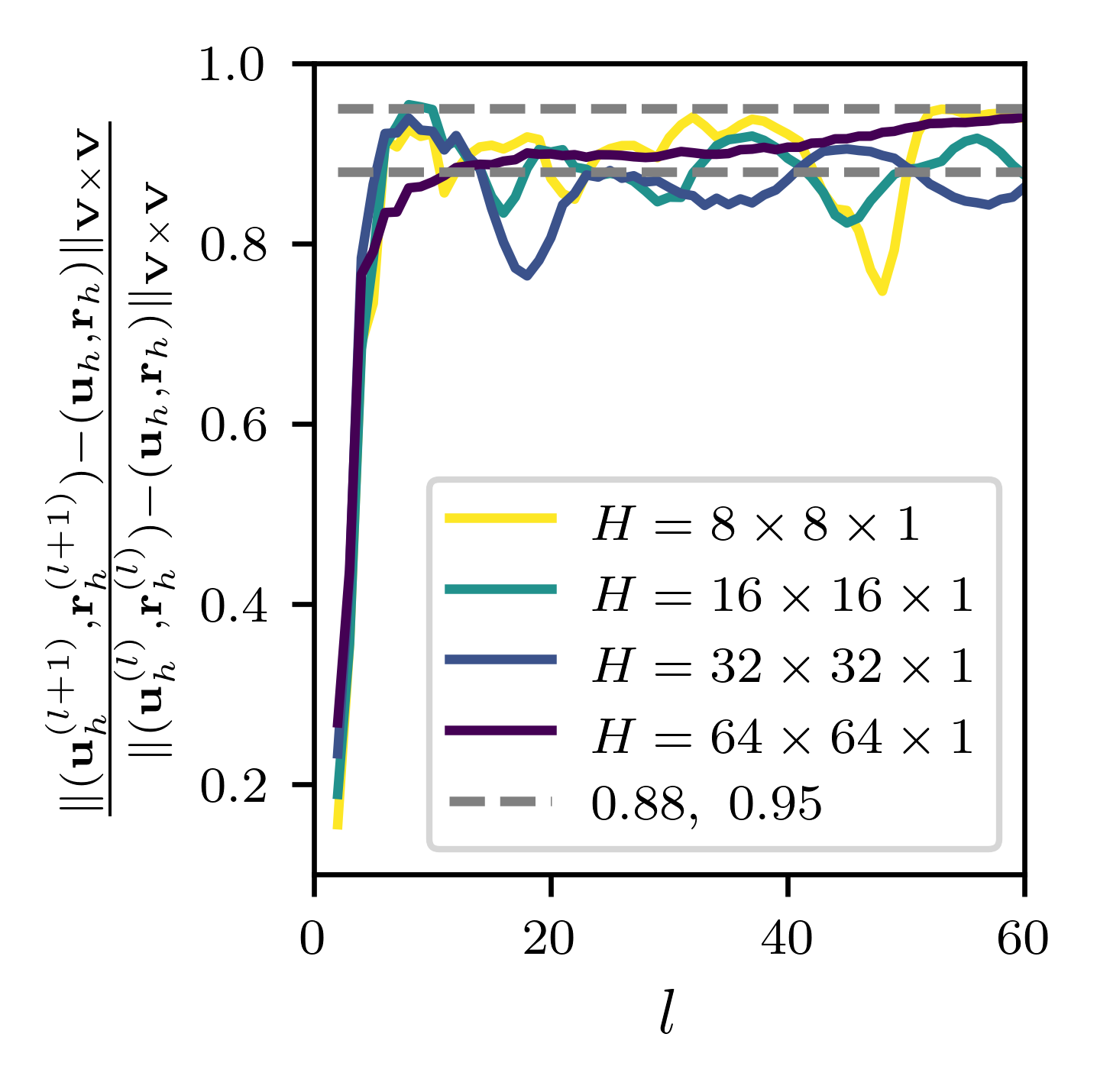}
\caption{Relative error as a function of the number of conjugate gradient iterations (left) and the error reduction factor versus the number of iterations (right).} 
\label{fig:fiber_h_convergence}
\end{figure}

The convergence rate predicted in \cref{thm:cg} is assessed by applying the domain decomposition method to the elliptic problem for several $H$. The initial guess $(\vec u_h^{(0)}, \vec r_h^{(0)}) \in \mathbf V_h \times \mathbf V_h$, consistent with the boundary conditions, is obtained by modeling the sheet as a continuum with straight Timoshenko beams, solving the continuum steady-state problem, and interpolating to the network. The results, presented in \cref{fig:fiber_h_convergence}, show linear convergence with no deterioration as $H$ decreases, confirming that \cref{thm:cg} applies for this example.

\subsubsection*{Domain decomposition with Crank--Nicolson}
\begin{figure}
\centering
\includegraphics[width = 0.8\textwidth]{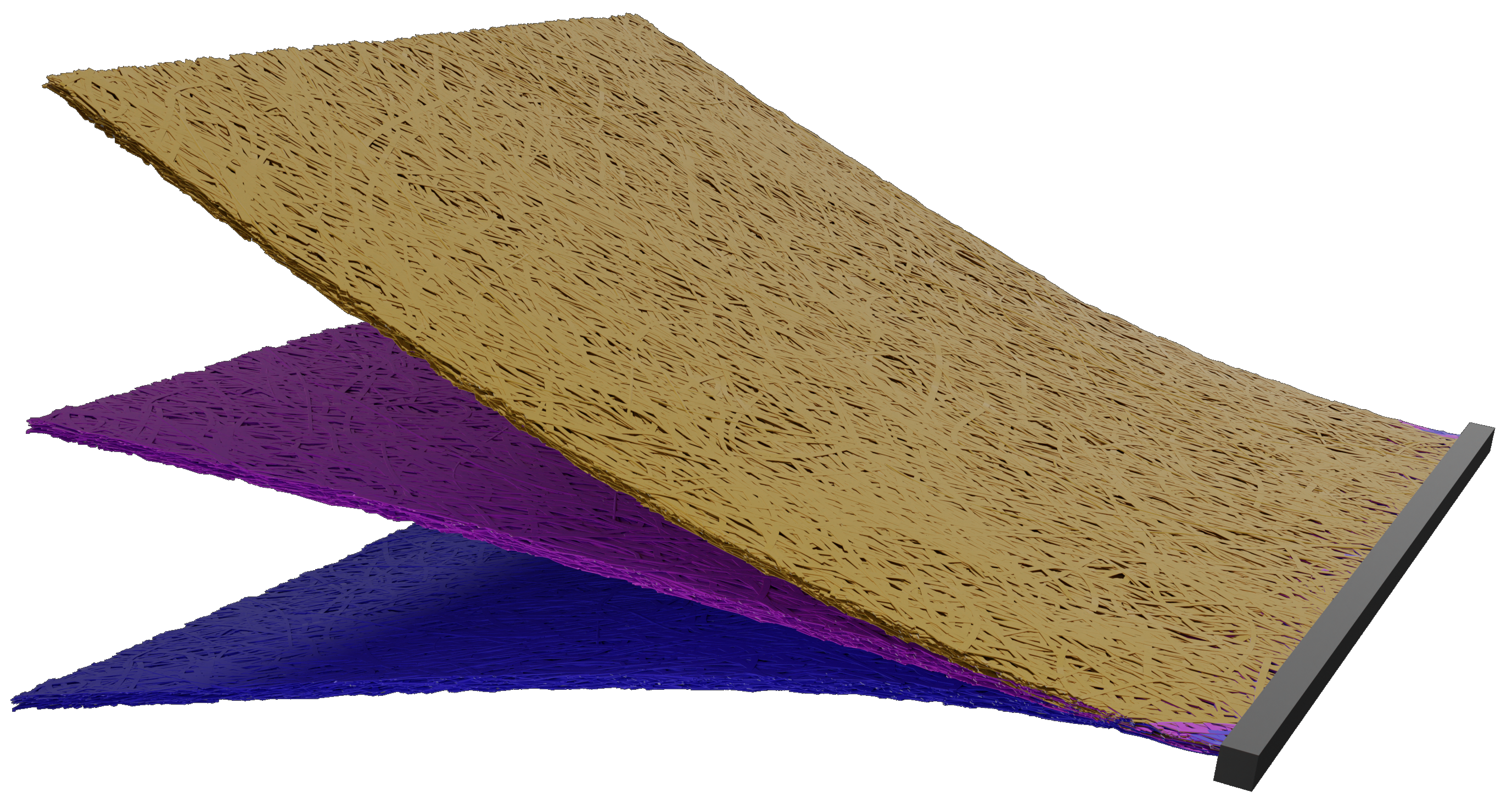}
\caption{The figure shows the solution at three time-steps in the simulation (yellow 0 ms, purple 0.242 ms, blue 0.484 ms).}
\label{fig:paper-wave}
\end{figure}

\begin{figure}
\centering
\includegraphics{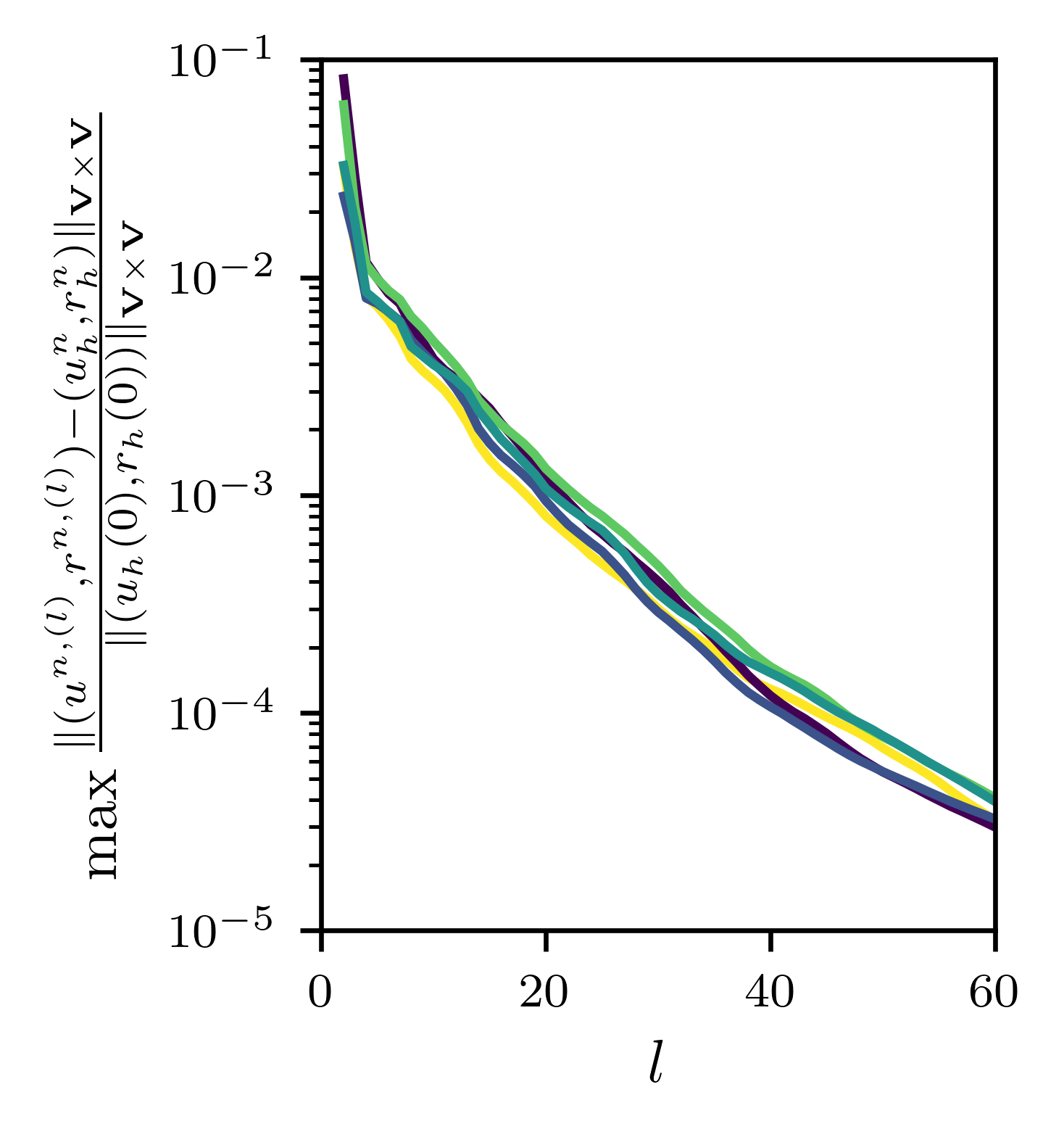}
\includegraphics{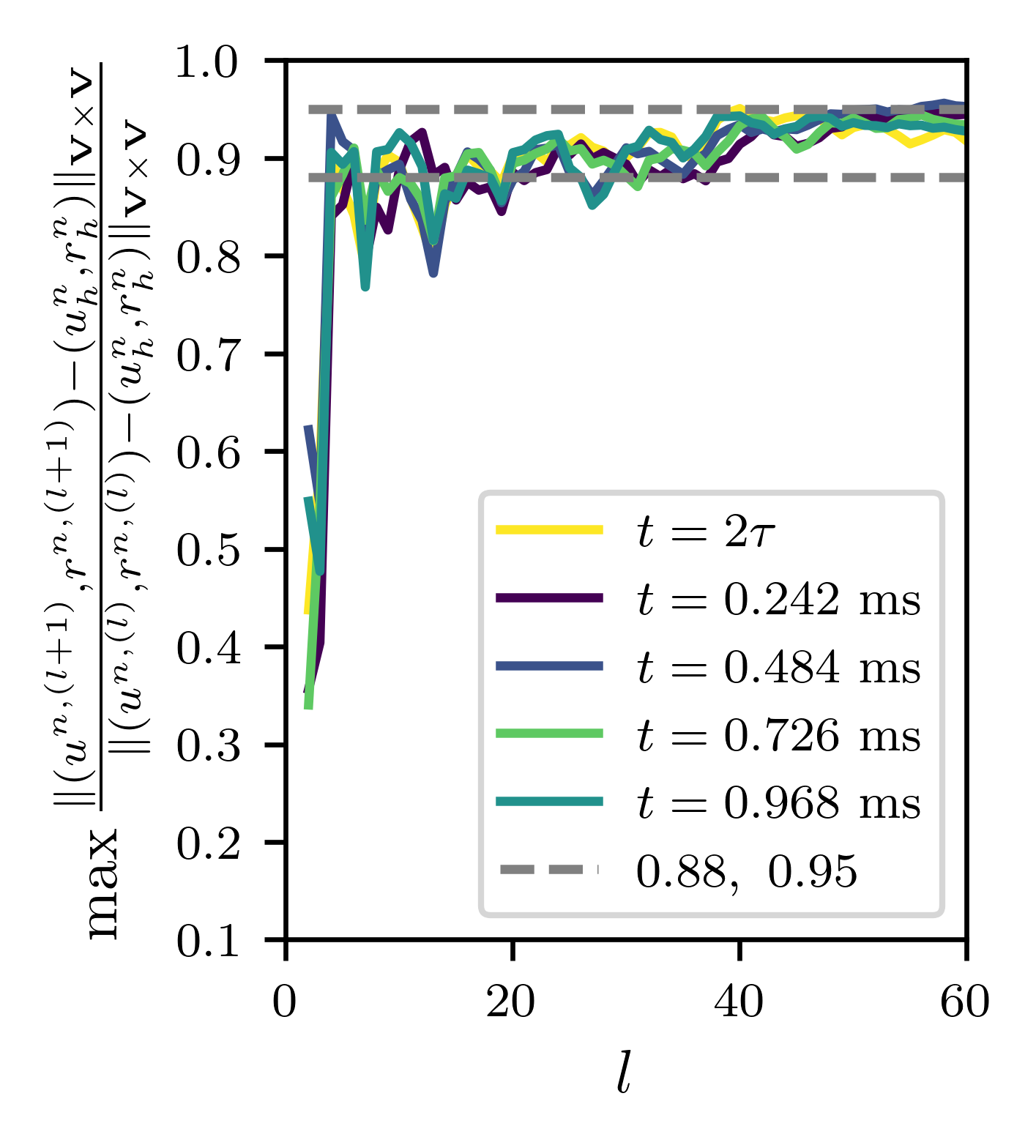}
\caption{The left figure shows the relative error at each iteration for the four time steps. The right figure shows the error reduction per domain decomposition iteration ($64 \times 64 \times 1$ elements) over the same four time steps.
} 
\label{fig:paper-wave-dd-conv}
\end{figure}

Next, we simulate the response of the paper after releasing the displacement Dirichlet condition by solving \eqref{eq:wave-tau-conv} using the Crank--Nicolson method. The material response closely resembles the first eigenmode, allowing an estimate of the period using data from \cref{tab:model-paper-prop}:
\begin{align*}
\text{Period}^{-1} &= \frac{(1.875)^2}{2\pi (8~\text{mm})^2} 
\sqrt{\frac{E_s I_x}{\rho_s (t_s \cdot (8~\text{mm}))}} \approx (1.1~\text{ms})^{-1}, \\
I_x &= \frac{t_s^3 (8~\text{mm})}{12}, \quad E_s = 12~\text{GPa},
\end{align*}
cf.~\cite[p.~484]{Gibson}. We set $T = 1.1~\text{ms}$ and simulate the response with time step $\tau = T/100$ using a direct linear solver. \cref{fig:paper-wave} shows the solution at three points in the first half of the period. The simulation is then repeated using the domain decomposition method with $64 \times 64 \times 1$ elements. The domain decomposition iterations are compared with the direct solver results, as shown in \cref{fig:paper-wave-dd-conv}. We again observe linear convergence of the preconditioned conjugate gradient method, essentially independent of the chosen time step.

\section{Conclusions and future work}\label{sec:conclusion}
We established optimal a priori error estimates for the proposed $\theta$-finite-element scheme applied to elastic wave propagation in Timoshenko beam network models with rigid joints. The theoretical results were validated through engineering numerical examples. For the explicit leapfrog method ($\theta = 0$), we also investigated the CFL condition both theoretically and numerically. Additionally, we proposed a subspace decomposition preconditioner for the efficient iterative solution of the arising linear systems ($\theta > 0$) and analyzed its convergence properties.
Future work will focus on multilevel preconditioners with algebraic coarsening and reduced overlap. A natural extension of the model is the inclusion of non-linear effects, incorporating geometrically non-linear formulations to account for large deformations, as well as non-linear constitutive laws to capture plasticity.

\section*{Acknowledgments}

M.~Hauck acknowledges funding from the Deutsche Forschungsgemeinschaft\linebreak  (DFG, German Research Foundation) -- Project-ID 258734477 -- SFB 1173. Furthermore,  A.~M\aa lqvist and L.~Swoboda acknowledge funding from the Swedish Research Council (VR) -- Project-ID 2023-03258\_VR.
\appendix

\bibliographystyle{ARalpha}
\bibliography{bib.bib}

\end{document}